
\input amstex

\expandafter\ifx\csname beta.def\endcsname\relax \else\endinput\fi
\expandafter\edef\csname beta.def\endcsname{%
 \catcode`\noexpand\@=\the\catcode`\@\space}

\let\atbefore @

\catcode`\@=11

\overfullrule\z@
\hsize 6.25truein
\vsize 9.63truein

\let\@ft@\expandafter \let\@tb@f@\atbefore

\newif\ifMag
\ifnum\mag>1000 \Magtrue\fi
=\ifMag cmr8\else cmr9\fi

\newdimen\p@@ \p@@\p@
\def\m@ths@r{\ifnum\mathsurround=\z@\z@\else\maths@r\fi}
\def\maths@r{1.6\p@@} \def\mathsurzero{\def\maths@r{\z@}}

\mathsurround\maths@r
\font\Brm=cmr12 \font\Bbf=cmbx12 \font\Bit=cmti12 \font\ssf=cmss10
\font\Bsl=cmsl10 scaled 1200 \font\Bmmi=cmmi10 scaled 1200
\font\BBf=cmbx12 scaled 1200 \font\BMmi=cmmi10 scaled 1440

\def\atletter{\edef\atrestore{\catcode`\noexpand\@=\the\catcode`\@\space}
 \catcode`\@=11}

\newread\@ux \newwrite\@@x \newwrite\@@cd
\let\@np@@\input
\def\@np@t#1{\openin\@ux#1\relax\ifeof\@ux\else\closein\@ux\relax\@np@@ #1\fi}
\def\input#1 {\openin\@ux#1\relax\ifeof\@ux\wrs@x{No file #1}\else
 \closein\@ux\relax\@np@@ #1\fi}
\def\Input#1 {\relax} 

\def\wr@@x#1{} \def\wrs@x{\immediate\write\sixt@@n}

\def\readldf{\@np@t{\jobname.ldf}}
\def\writeldf{\def\wr@@x{\immediate\write\@@x}
 \def\cl@selbl{\wr@@x{\string\Snodef{\the\Sno}}\wr@@x{\string\endinput}%
 \immediate\closeout\@@x} \immediate\openout\@@x\jobname.ldf}
\let\cl@selbl\relax

\def\tod@y{\ifcase\month\or
 January\or February\or March\or April\or May\or June\or July\or
 August\or September\or October\or November\or December\fi\space\,
\number\day,\space\,\number\year}

\newcount\c@time
\def\h@@r{hh}\def\m@n@te{mm}
\def\wh@tt@me{\c@time\time\divide\c@time 60\edef\h@@r{\number\c@time}%
 \multiply\c@time -60\advance\c@time\time\edef
 \m@n@te{\ifnum\c@time<10 0\fi\number\c@time}}
\def\t@me{\h@@r\/{\rm:}\m@n@te}  \let\whattime\wh@tt@me
\def\today{\tod@y\wr@@x{\string\todaydef{\tod@y}}}
\def\nowtime{\t@me{\let\/\ic@\wr@@x{\string\nowtimedef{\t@me}}}}
\def\todaydef#1{} \def\nowtimedef#1{}

\def\em#1{{\it #1\/}} \def\emph#1{{\sl #1\/}}

\def\fitem#1{\par\setbox\z@\hbox{#1}\hangindent\wd\z@
 \hglue-2\parindent\kern\wd\z@\indent\llap{#1}\ignore}

\def\itemflat#1{\par\setbox\z@\hbox{\rm #1\enspace}\hang\ifnum\wd\z@>\parindent
 \noindent\unhbox\z@\ignore\else\textindent{\rm#1}\fi}

\newcount\itemlet
\def\newbi{\itemlet 96} \newbi
\def\bitem{\gad\itemlet \par\hangindent1.5\parindent
 \hglue-.5\parindent\textindent{\rm\rlap{\char\the\itemlet}\hp{b})}}

\newcount\itemrm

\def\iitem{\gad\itemrm \par\hangindent1.5\parindent
 \hglue-.5\parindent\textindent{\rm\hp{v}\llap{\romannumeral\the\itemrm})}}

\def\center{\par\begingroup\leftskip\z@ plus \hsize \rightskip\leftskip
 \parindent\z@\parfillskip\z@skip \def\\{\unskip\break}}
\def\endcenter{\endgraf\endgroup}

\def\Abstract{\begingroup\narrower\nt{\bf Abstract.}\enspace\ignore}
\def\endAbs{\endgraf\endgroup}

\let\b@gr@@\begingroup \let\B@gr@@\begingroup
\def\b@gr@{\b@gr@@\let\b@gr@@\undefined}
\def\B@gr@{\B@gr@@\let\B@gr@@\undefined}

\def\@fn@xt#1#2#3{\let\@ch@r=#1\def\n@xt{\ifx\t@st@\@ch@r
 \def\n@@xt{#2}\else\def\n@@xt{#3}\fi\n@@xt}\futurelet\t@st@\n@xt}

\def\@fwd@@#1#2#3{\setbox\z@\hbox{#1}\ifdim\wd\z@>\z@#2\else#3\fi}
\def\s@twd@#1#2{\setbox\z@\hbox{#2}#1\wd\z@}

\def\r@st@re#1{\let#1\s@v@} \def\s@v@d@f{\let\s@v@}

\def\p@sk@p#1#2{\par\skip@#2\relax\ifdim\lastskip<\skip@\relax\removelastskip
 \ifnum#1=\z@\else\penalty#1\relax\fi\vskip\skip@
 \else\ifnum#1=\z@\else\penalty#1\relax\fi\fi}
\def\sk@@p#1{\par\skip@#1\relax\ifdim\lastskip<\skip@\relax\removelastskip
 \vskip\skip@\fi}

\newbox\p@b@ld
\def\poorbold#1{\setbox\p@b@ld\hbox{#1}\kern-.01em\copy\p@b@ld\kern-\wd\p@b@ld
 \kern.02em\copy\p@b@ld\kern-\wd\p@b@ld\kern-.012em\raise.02em\box\p@b@ld}

\ifx\plainfootnote\undefined \let\plainfootnote\footnote \fi

\let\s@v@\proclaim \let\proclaim\relax
\def\r@R@fs#1{\let#1\s@R@fs} \let\s@R@fs\Refs \let\Refs\relax
\def\r@endd@#1{\let#1\s@endd@} \let\s@endd@\enddocument
\let\bye\relax

\def\myR@fs{\@fn@xt[\m@R@f@\m@R@fs} \def\m@R@fs{\@fn@xt*\m@r@f@@\m@R@f@@}
\def\m@R@f@@{\m@R@f@[References]} \def\m@r@f@@*{\m@R@f@[]}

\def\Twelvepoint{\twelvepoint \let\Bbf\BBf \let\Bmmi\BMmi
\font\Brm=cmr12 scaled 1200 \font\Bit=cmti12 scaled 1200
\font\ssf=cmss10 scaled 1200 \font\Bsl=cmsl10 scaled 1440
\font\BBf=cmbx12 scaled 1440 \font\BMmi=cmmi10 scaled 1728}

\newdimen\b@gsize

\newdimen\r@f@nd \newbox\r@f@b@x \newbox\adjb@x
\newbox\p@nct@ \newbox\k@yb@x \newcount\rcount
\newbox\b@b@x \newbox\p@p@rb@x \newbox\j@@rb@x \newbox\y@@rb@x
\newbox\v@lb@x \newbox\is@b@x \newbox\p@g@b@x \newif\ifp@g@ \newif\ifp@g@s
\newbox\inb@@kb@x \newbox\b@@kb@x \newbox\p@blb@x \newbox\p@bl@db@x
\newbox\ed@b@x \newif\ifed@ \newif\ifed@s \newif\if@fl@b \newif\if@fn@m
\newbox\p@p@nf@b@x \newbox\inf@b@x \newbox\b@@nf@b@x

\newif\ifp@gen@

\newif\ifamsppt

\@ft@\ifx\csname amsppt.sty\endcsname\relax

\headline={\hfil}
\footline={\ifp@gen@\ifnum\pageno=\z@\else\hfil\foliorm\folio\fi\else
 \ifnum\pageno=\z@\hfil\foliorm\folio\fi\fi\hfil\global\p@gen@true}
\parindent1pc

\font@\tensmc=cmcsc10
\font@\sevenex=cmex7
\font@\sevenit=cmti7
\font@\eightrm=cmr8
\font@\sixrm=cmr6
\font@\eighti=cmmi8 \skewchar\eighti='177
\font@\sixi=cmmi6 \skewchar\sixi='177
\font@\eightsy=cmsy8 \skewchar\eightsy='60
\font@\sixsy=cmsy6 \skewchar\sixsy='60
\font@\eightex=cmex8
\font@\eightbf=cmbx8
\font@\sixbf=cmbx6
\font@\eightit=cmti8
\font@\eightsl=cmsl8
\font@\eightsmc=cmcsc8
\font@\eighttt=cmtt8
\font@\ninerm=cmr9
\font@\ninei=cmmi9 \skewchar\ninei='177
\font@\ninesy=cmsy9 \skewchar\ninesy='60
\font@\nineex=cmex9
\font@\ninebf=cmbx9
\font@\nineit=cmti9
\font@\ninesl=cmsl9
\font@\ninesmc=cmcsc9
\font@\ninemsa=msam9
\font@\ninemsb=msbm9
\font@\nineeufm=eufm9
\font@\eightmsa=msam8
\font@\eightmsb=msbm8
\font@\eighteufm=eufm8
\font@\sixmsa=msam6
\font@\sixmsb=msbm6
\font@\sixeufm=eufm6

\loadmsam\loadmsbm\loadeufm
\input amssym.tex

\def\footnoterule{\kern-3\p@\hrule width5pc\kern 2.6\p@}
\def\m@k@foot#1{\insert\footins
 {\interlinepenalty\interfootnotelinepenalty
 \eightpoint\splittopskip\ht\strutbox\splitmaxdepth\dp\strutbox
 \floatingpenalty\@MM\leftskip\z@\rightskip\z@
 \spaceskip\z@\xspaceskip\z@
 \leavevmode\footstrut\ignore#1\unskip\lower\dp\strutbox
 \vbox to\dp\strutbox{}}}
\def\ftext#1{\m@k@foot{\vsk-.8>\nt #1}}
\def\pr@cl@@m#1{\p@sk@p{-100}\medskipamount\b@gr@\nt\ignore
 \bf #1\unskip.\enspace\sl\ignore}
\outer\def\proclaim{\pr@cl@@m} \s@v@d@f\proclaim \let\proclaim\relax
\def\endproclaim{\endgroup\p@sk@p{55}\medskipamount}
\def\demo#1{\sk@@p\medskipamount\nt{\ignore\it #1\unskip.}\enspace
 \ignore}
\def\enddemo{\sk@@p\medskipamount}

\def\cite#1{{\rm[#1]}} \let\nofrills\relax
 \def\Refs#1#2{\relax}

\def\big@#1#2{{\hbox{$\left#2\vcenter to#1\b@gsize{}%
 \right.\nulldelimiterspace\z@\m@th$}}}
\def\big{\big@\@ne}
\def\Big{\big@{1.5}}
\def\bigg{\big@\tw@}
\def\Bigg{\big@{2.5}}
\normallineskiplimit\p@

\def\tenpoint{\p@@\p@ \normallineskiplimit\p@@
 \mathsurround\m@ths@r \normalbaselineskip12\p@@
 \abovedisplayskip12\p@@ plus3\p@@ minus9\p@@
 \belowdisplayskip\abovedisplayskip
 \abovedisplayshortskip\z@ plus3\p@@
 \belowdisplayshortskip7\p@@ plus3\p@@ minus4\p@@
 \textonlyfont@\rm\tenrm \textonlyfont@\it\tenit
 \textonlyfont@\sl\tensl \textonlyfont@\bf\tenbf
 \textonlyfont@\smc\tensmc \textonlyfont@\tt\tentt
 \ifsyntax@ \def\big##1{{\hbox{$\left##1\right.$}}}%
  \let\Big\big \let\bigg\big \let\Bigg\big
 \else
  \textfont\z@\tenrm \scriptfont\z@\sevenrm \scriptscriptfont\z@\fiverm
  \textfont\@ne\teni \scriptfont\@ne\seveni \scriptscriptfont\@ne\fivei
  \textfont\tw@\tensy \scriptfont\tw@\sevensy \scriptscriptfont\tw@\fivesy
  \textfont\thr@@\tenex \scriptfont\thr@@\sevenex
	\scriptscriptfont\thr@@\sevenex
  \textfont\itfam\tenit \scriptfont\itfam\sevenit
	\scriptscriptfont\itfam\sevenit
  \textfont\bffam\tenbf \scriptfont\bffam\sevenbf
	\scriptscriptfont\bffam\fivebf
  \textfont\msafam\tenmsa \scriptfont\msafam\sevenmsa
	\scriptscriptfont\msafam\fivemsa
  \textfont\msbfam\tenmsb \scriptfont\msbfam\sevenmsb
	\scriptscriptfont\msbfam\fivemsb
  \textfont\eufmfam\teneufm \scriptfont\eufmfam\seveneufm
	\scriptscriptfont\eufmfam\fiveeufm
  \setbox\strutbox\hbox{\vrule height8.5\p@@ depth3.5\p@@ width\z@}%
  \setbox\strutbox@\hbox{\lower.5\normallineskiplimit\vbox{%
	\kern-\normallineskiplimit\copy\strutbox}}%
   \setbox\z@\vbox{\hbox{$($}\kern\z@}\b@gsize1.2\ht\z@
  \fi
  \normalbaselines\rm\dotsspace@1.5mu\ex@.2326ex\jot3\ex@}

\def\eightpoint{\p@@.8\p@ \normallineskiplimit\p@@
 \mathsurround\m@ths@r \normalbaselineskip10\p@
 \abovedisplayskip10\p@ plus2.4\p@ minus7.2\p@
 \belowdisplayskip\abovedisplayskip
 \abovedisplayshortskip\z@ plus3\p@@
 \belowdisplayshortskip7\p@@ plus3\p@@ minus4\p@@
 \textonlyfont@\rm\eightrm \textonlyfont@\it\eightit
 \textonlyfont@\sl\eightsl \textonlyfont@\bf\eightbf
 \textonlyfont@\smc\eightsmc \textonlyfont@\tt\eighttt
 \ifsyntax@\def\big##1{{\hbox{$\left##1\right.$}}}%
  \let\Big\big \let\bigg\big \let\Bigg\big
 \else
  \textfont\z@\eightrm \scriptfont\z@\sixrm \scriptscriptfont\z@\fiverm
  \textfont\@ne\eighti \scriptfont\@ne\sixi \scriptscriptfont\@ne\fivei
  \textfont\tw@\eightsy \scriptfont\tw@\sixsy \scriptscriptfont\tw@\fivesy
  \textfont\thr@@\eightex \scriptfont\thr@@\sevenex
	\scriptscriptfont\thr@@\sevenex
  \textfont\itfam\eightit \scriptfont\itfam\sevenit
	\scriptscriptfont\itfam\sevenit
  \textfont\bffam\eightbf \scriptfont\bffam\sixbf
	\scriptscriptfont\bffam\fivebf
  \textfont\msafam\eightmsa \scriptfont\msafam\sixmsa
	\scriptscriptfont\msafam\fivemsa
  \textfont\msbfam\eightmsb \scriptfont\msbfam\sixmsb
	\scriptscriptfont\msbfam\fivemsb
  \textfont\eufmfam\eighteufm \scriptfont\eufmfam\sixeufm
	\scriptscriptfont\eufmfam\fiveeufm
 \setbox\strutbox\hbox{\vrule height7\p@ depth3\p@ width\z@}%
 \setbox\strutbox@\hbox{\raise.5\normallineskiplimit\vbox{%
   \kern-\normallineskiplimit\copy\strutbox}}%
 \setbox\z@\vbox{\hbox{$($}\kern\z@}\b@gsize1.2\ht\z@
 \fi
 \normalbaselines\eightrm\dotsspace@1.5mu\ex@.2326ex\jot3\ex@}

\def\ninepoint{\p@@.9\p@ \normallineskiplimit\p@@
 \mathsurround\m@ths@r \normalbaselineskip11\p@
 \abovedisplayskip11\p@ plus2.7\p@ minus8.1\p@
 \belowdisplayskip\abovedisplayskip
 \abovedisplayshortskip\z@ plus3\p@@
 \belowdisplayshortskip7\p@@ plus3\p@@ minus4\p@@
 \textonlyfont@\rm\ninerm \textonlyfont@\it\nineit
 \textonlyfont@\sl\ninesl \textonlyfont@\bf\ninebf
 \textonlyfont@\smc\ninesmc \textonlyfont@\tt\ninett
 \ifsyntax@ \def\big##1{{\hbox{$\left##1\right.$}}}%
  \let\Big\big \let\bigg\big \let\Bigg\big
 \else
  \textfont\z@\ninerm \scriptfont\z@\sevenrm \scriptscriptfont\z@\fiverm
  \textfont\@ne\ninei \scriptfont\@ne\seveni \scriptscriptfont\@ne\fivei
  \textfont\tw@\ninesy \scriptfont\tw@\sevensy \scriptscriptfont\tw@\fivesy
  \textfont\thr@@\nineex \scriptfont\thr@@\sevenex
	\scriptscriptfont\thr@@\sevenex
  \textfont\itfam\nineit \scriptfont\itfam\sevenit
	\scriptscriptfont\itfam\sevenit
  \textfont\bffam\ninebf \scriptfont\bffam\sevenbf
	\scriptscriptfont\bffam\fivebf
  \textfont\msafam\ninemsa \scriptfont\msafam\sevenmsa
	\scriptscriptfont\msafam\fivemsa
  \textfont\msbfam\ninemsb \scriptfont\msbfam\sevenmsb
	\scriptscriptfont\msbfam\fivemsb
  \textfont\eufmfam\nineeufm \scriptfont\eufmfam\seveneufm
	\scriptscriptfont\eufmfam\fiveeufm
  \setbox\strutbox\hbox{\vrule height8.5\p@@ depth3.5\p@@ width\z@}%
  \setbox\strutbox@\hbox{\lower.5\normallineskiplimit\vbox{%
	\kern-\normallineskiplimit\copy\strutbox}}%
   \setbox\z@\vbox{\hbox{$($}\kern\z@}\b@gsize1.2\ht\z@
  \fi
  \normalbaselines\rm\dotsspace@1.5mu\ex@.2326ex\jot3\ex@}

\font@\twelverm=cmr10 scaled 1200
\font@\twelveit=cmti10 scaled 1200
\font@\twelvesl=cmsl10 scaled 1200
\font@\twelvebf=cmbx10 scaled 1200
\font@\twelvesmc=cmcsc10 scaled 1200
\font@\twelvett=cmtt10 scaled 1200
\font@\twelvei=cmmi10 scaled 1200 \skewchar\twelvei='177
\font@\twelvesy=cmsy10 scaled 1200 \skewchar\twelvesy='60
\font@\twelveex=cmex10 scaled 1200
\font@\twelvemsa=msam10 scaled 1200
\font@\twelvemsb=msbm10 scaled 1200
\font@\twelveeufm=eufm10 scaled 1200

\def\twelvepoint{\p@@1.2\p@ \normallineskiplimit\p@@
 \mathsurround\m@ths@r \normalbaselineskip12\p@@
 \abovedisplayskip12\p@@ plus3\p@@ minus9\p@@
 \belowdisplayskip\abovedisplayskip
 \abovedisplayshortskip\z@ plus3\p@@
 \belowdisplayshortskip7\p@@ plus3\p@@ minus4\p@@
 \textonlyfont@\rm\twelverm \textonlyfont@\it\twelveit
 \textonlyfont@\sl\twelvesl \textonlyfont@\bf\twelvebf
 \textonlyfont@\smc\twelvesmc \textonlyfont@\tt\twelvett
 \ifsyntax@ \def\big##1{{\hbox{$\left##1\right.$}}}%
  \let\Big\big \let\bigg\big \let\Bigg\big
 \else
  \textfont\z@\twelverm \scriptfont\z@\eightrm \scriptscriptfont\z@\sixrm
  \textfont\@ne\twelvei \scriptfont\@ne\eighti \scriptscriptfont\@ne\sixi
  \textfont\tw@\twelvesy \scriptfont\tw@\eightsy \scriptscriptfont\tw@\sixsy
  \textfont\thr@@\twelveex \scriptfont\thr@@\eightex
	\scriptscriptfont\thr@@\sevenex
  \textfont\itfam\twelveit \scriptfont\itfam\eightit
	\scriptscriptfont\itfam\sevenit
  \textfont\bffam\twelvebf \scriptfont\bffam\eightbf
	\scriptscriptfont\bffam\sixbf
  \textfont\msafam\twelvemsa \scriptfont\msafam\eightmsa
	\scriptscriptfont\msafam\sixmsa
  \textfont\msbfam\twelvemsb \scriptfont\msbfam\eightmsb
	\scriptscriptfont\msbfam\sixmsb
  \textfont\eufmfam\twelveeufm \scriptfont\eufmfam\eighteufm
	\scriptscriptfont\eufmfam\sixeufm
  \setbox\strutbox\hbox{\vrule height8.5\p@@ depth3.5\p@@ width\z@}%
  \setbox\strutbox@\hbox{\lower.5\normallineskiplimit\vbox{%
	\kern-\normallineskiplimit\copy\strutbox}}%
  \setbox\z@\vbox{\hbox{$($}\kern\z@}\b@gsize1.2\ht\z@
  \fi
  \normalbaselines\rm\dotsspace@1.5mu\ex@.2326ex\jot3\ex@}

\font@\twelvetrm=cmr10 at 12truept
\font@\twelvetit=cmti10 at 12truept
\font@\twelvetsl=cmsl10 at 12truept
\font@\twelvetbf=cmbx10 at 12truept
\font@\twelvetsmc=cmcsc10 at 12truept
\font@\twelvettt=cmtt10 at 12truept
\font@\twelveti=cmmi10 at 12truept \skewchar\twelveti='177
\font@\twelvetsy=cmsy10 at 12truept \skewchar\twelvetsy='60
\font@\twelvetex=cmex10 at 12truept
\font@\twelvetmsa=msam10 at 12truept
\font@\twelvetmsb=msbm10 at 12truept
\font@\twelveteufm=eufm10 at 12truept

\def\twelvetruepoint{\p@@1.2truept \normallineskiplimit\p@@
 \mathsurround\m@ths@r \normalbaselineskip12\p@@
 \abovedisplayskip12\p@@ plus3\p@@ minus9\p@@
 \belowdisplayskip\abovedisplayskip
 \abovedisplayshortskip\z@ plus3\p@@
 \belowdisplayshortskip7\p@@ plus3\p@@ minus4\p@@
 \textonlyfont@\rm\twelvetrm \textonlyfont@\it\twelvetit
 \textonlyfont@\sl\twelvetsl \textonlyfont@\bf\twelvetbf
 \textonlyfont@\smc\twelvetsmc \textonlyfont@\tt\twelvettt
 \ifsyntax@ \def\big##1{{\hbox{$\left##1\right.$}}}%
  \let\Big\big \let\bigg\big \let\Bigg\big
 \else
  \textfont\z@\twelvetrm \scriptfont\z@\eightrm \scriptscriptfont\z@\sixrm
  \textfont\@ne\twelveti \scriptfont\@ne\eighti \scriptscriptfont\@ne\sixi
  \textfont\tw@\twelvetsy \scriptfont\tw@\eightsy \scriptscriptfont\tw@\sixsy
  \textfont\thr@@\twelvetex \scriptfont\thr@@\eightex
	\scriptscriptfont\thr@@\sevenex
  \textfont\itfam\twelvetit \scriptfont\itfam\eightit
	\scriptscriptfont\itfam\sevenit
  \textfont\bffam\twelvetbf \scriptfont\bffam\eightbf
	\scriptscriptfont\bffam\sixbf
  \textfont\msafam\twelvetmsa \scriptfont\msafam\eightmsa
	\scriptscriptfont\msafam\sixmsa
  \textfont\msbfam\twelvetmsb \scriptfont\msbfam\eightmsb
	\scriptscriptfont\msbfam\sixmsb
  \textfont\eufmfam\twelveteufm \scriptfont\eufmfam\eighteufm
	\scriptscriptfont\eufmfam\sixeufm
  \setbox\strutbox\hbox{\vrule height8.5\p@@ depth3.5\p@@ width\z@}%
  \setbox\strutbox@\hbox{\lower.5\normallineskiplimit\vbox{%
	\kern-\normallineskiplimit\copy\strutbox}}%
  \setbox\z@\vbox{\hbox{$($}\kern\z@}\b@gsize1.2\ht\z@
  \fi
  \normalbaselines\rm\dotsspace@1.5mu\ex@.2326ex\jot3\ex@}

\font@\elevenrm=cmr10 scaled 1095
\font@\elevenit=cmti10 scaled 1095
\font@\elevensl=cmsl10 scaled 1095
\font@\elevenbf=cmbx10 scaled 1095
\font@\elevensmc=cmcsc10 scaled 1095
\font@\eleventt=cmtt10 scaled 1095
\font@\eleveni=cmmi10 scaled 1095 \skewchar\eleveni='177
\font@\elevensy=cmsy10 scaled 1095 \skewchar\elevensy='60
\font@\elevenex=cmex10 scaled 1095
\font@\elevenmsa=msam10 scaled 1095
\font@\elevenmsb=msbm10 scaled 1095
\font@\eleveneufm=eufm10 scaled 1095

\def\elevenpoint{\p@@1.1\p@ \normallineskiplimit\p@@
 \mathsurround\m@ths@r \normalbaselineskip12\p@@
 \abovedisplayskip12\p@@ plus3\p@@ minus9\p@@
 \belowdisplayskip\abovedisplayskip
 \abovedisplayshortskip\z@ plus3\p@@
 \belowdisplayshortskip7\p@@ plus3\p@@ minus4\p@@
 \textonlyfont@\rm\elevenrm \textonlyfont@\it\elevenit
 \textonlyfont@\sl\elevensl \textonlyfont@\bf\elevenbf
 \textonlyfont@\smc\elevensmc \textonlyfont@\tt\eleventt
 \ifsyntax@ \def\big##1{{\hbox{$\left##1\right.$}}}%
  \let\Big\big \let\bigg\big \let\Bigg\big
 \else
  \textfont\z@\elevenrm \scriptfont\z@\eightrm \scriptscriptfont\z@\sixrm
  \textfont\@ne\eleveni \scriptfont\@ne\eighti \scriptscriptfont\@ne\sixi
  \textfont\tw@\elevensy \scriptfont\tw@\eightsy \scriptscriptfont\tw@\sixsy
  \textfont\thr@@\elevenex \scriptfont\thr@@\eightex
	\scriptscriptfont\thr@@\sevenex
  \textfont\itfam\elevenit \scriptfont\itfam\eightit
	\scriptscriptfont\itfam\sevenit
  \textfont\bffam\elevenbf \scriptfont\bffam\eightbf
	\scriptscriptfont\bffam\sixbf
  \textfont\msafam\elevenmsa \scriptfont\msafam\eightmsa
	\scriptscriptfont\msafam\sixmsa
  \textfont\msbfam\elevenmsb \scriptfont\msbfam\eightmsb
	\scriptscriptfont\msbfam\sixmsb
  \textfont\eufmfam\eleveneufm \scriptfont\eufmfam\eighteufm
	\scriptscriptfont\eufmfam\sixeufm
  \setbox\strutbox\hbox{\vrule height8.5\p@@ depth3.5\p@@ width\z@}%
  \setbox\strutbox@\hbox{\lower.5\normallineskiplimit\vbox{%
	\kern-\normallineskiplimit\copy\strutbox}}%
  \setbox\z@\vbox{\hbox{$($}\kern\z@}\b@gsize1.2\ht\z@
  \fi
  \normalbaselines\rm\dotsspace@1.5mu\ex@.2326ex\jot3\ex@}

\def\m@R@f@[#1]{\mathsurzero{
 \s@ct{}{#1}}\wr@@c{\string\Refcd{#1}{\the\pageno}}\B@gr@
 \frenchspacing\rcount\z@\refkey{[##1]}\refno{[##1]}\widest{AZ}\keyright
 \let\Key\key\let\refin\relax}
\def\widest#1{\s@twd@\r@f@nd{\r@fk@y{#1}\enspace}}
\def\widestno#1{\s@twd@\r@f@nd{\r@fn@{#1}\enspace}}
\def\widestlabel#1{\s@twd@\r@f@nd{#1\enspace}}
\def\refkey{\def\r@fk@y##1} \def\refno{\def\r@fn@##1}
\def\keyright{\def\r@fit@m{\hang\textindent}}
\def\keyflat{\def\r@fit@m##1{\setbox\z@\hbox{\rm ##1\enspace}\hang\noindent
 \ifnum\wd\z@<\parindent\indent\hglue-\wd\z@\fi\unhbox\z@}}

\def\R@fb@x{\global\setbox\r@f@b@x} \def\K@yb@x{\global\setbox\k@yb@x}
\def\ref{\par\b@gr@\rm\R@fb@x\box\voidb@x\K@yb@x\box\voidb@x\@fn@mfalse
 \@fl@bfalse\b@g@nr@f}
\def\c@nc@t#1{\setbox\z@\lastbox
 \setbox\adjb@x\hbox{\unhbox\adjb@x\unhbox\z@\unskip\unskip\unpenalty#1}}
\def\adjust#1{\relax\ifmmode\penalty-\@M\null\hfil$\clubpenalty\z@
 \widowpenalty\z@\interlinepenalty\z@\offinterlineskip\endgraf
 \setbox\z@\lastbox\unskip\unpenalty\c@nc@t{#1}\nt$\hfil\penalty-\@M
 \else\endgraf\c@nc@t{#1}\nt\fi}
\def\adjustnext#1{\P@nct\hbox{#1}\ignore}
\def\adjustend#1{\def\@djp@{#1}\ignore}

\def\cl@s@{\adjust{\@djp@}\endgraf\setbox\z@\lastbox
 \global\setbox\@ne\hbox{\unhbox\adjb@x\ifvoid\z@\else\unhbox\z@\unskip\unskip
 \unpenalty\fi}\egroup\ifnum\c@rr@nt=\k@yb@x\global\fi
 \setbox\c@rr@nt\hbox{\unhbox\@ne\box\p@nct@}\P@nct\null}
\def\@p@n#1{\def\c@rr@nt{#1}\setbox\c@rr@nt\vbox\bgroup\let\@djp@\relax
 \hsize\maxdimen\nt}
\def\b@g@nr@f{\bgroup\@p@n\z@}
\def\key{\cl@s@\ifvoid\k@yb@x\@p@n\k@yb@x\else\@p@n\z@\fi}
\def\label{\cl@s@\ifvoid\k@yb@x\global\@fl@btrue\@p@n\k@yb@x\else\@p@n\z@\fi}
\def\no{\cl@s@\ifvoid\k@yb@x\gad\rcount\global\@fn@mtrue
 \K@yb@x\hbox{\the\rcount}\fi\@p@n\z@}
\def\labelno{\cl@s@\ifvoid\k@yb@x\gad\rcount\@fl@btrue\@p@n\k@yb@x\the\rcount
 \else\@p@n\z@\fi}
\def\by{\cl@s@\@p@n\b@b@x} \def\paper{\cl@s@\@p@n\p@p@rb@x\it\ignore}
\def\jour{\cl@s@\@p@n\j@@rb@x} \def\yr{\cl@s@\@p@n\y@@rb@x}
\def\vol{\cl@s@\@p@n\v@lb@x\bf\ignore} \def\issue{\cl@s@\@p@n\is@b@x}
\def\page{\cl@s@\ifp@g@s\@p@n\z@\else\p@g@true\@p@n\p@g@b@x\fi}
\def\pages{\cl@s@\ifp@g@\@p@n\z@\else\p@g@strue\@p@n\p@g@b@x\fi}
\def\inbook{\cl@s@\@p@n\inb@@kb@x} \def\book{\cl@s@\@p@n\b@@kb@x\it\ignore}
\def\publ{\cl@s@\@p@n\p@blb@x} \def\publaddr{\cl@s@\@p@n\p@bl@db@x}
\def\ed{\cl@s@\ifed@s\@p@n\z@\else\ed@true\@p@n\ed@b@x\fi}
\def\eds{\cl@s@\ifed@\@p@n\z@\else\ed@strue\@p@n\ed@b@x\fi}
\def\info{\cl@s@\@p@n\inf@b@x} \def\paperinfo{\cl@s@\@p@n\p@p@nf@b@x}
\def\bookinfo{\cl@s@\@p@n\b@@nf@b@x} \let\finalinfo\info
\def\P@nct{\global\setbox\p@nct@} \def\nopunct{\P@nct\box\voidb@x}
\def\p@@@t#1#2{\ifvoid\p@nct@\else#1\unhbox\p@nct@#2\fi}
\def\sp@@{\penalty-50 \space\hskip\z@ plus.1em}
\def\c@mm@{\p@@@t,\sp@@} \def\sp@c@{\p@@@t\empty\sp@@} \def\p@@nt{.\kern.3em}
\def\p@tb@x#1#2{\ifvoid#1\else#2\@nb@x#1\fi}
\def\@nb@x#1{\unhbox#1\P@nct\lastbox}
\def\endr@f@{\cl@s@\nopunct
 \R@fb@x\hbox{\unhbox\r@f@b@x \p@tb@x\b@b@x\empty
 \ifvoid\j@@rb@x\ifvoid\inb@@kb@x\ifvoid\p@p@rb@x\ifvoid\b@@kb@x
  \ifvoid\p@p@nf@b@x\ifvoid\b@@nf@b@x
  \p@tb@x\v@lb@x\c@mm@ \ifvoid\y@@rb@x\else\sp@c@(\@nb@x\y@@rb@x)\fi
  \p@tb@x\is@b@x{\c@mm@ no\p@@nt}\p@tb@x\p@g@b@x\c@mm@ \p@tb@x\inf@b@x\c@mm@
  \else\p@tb@x \b@@nf@b@x\c@mm@ \p@tb@x\v@lb@x\c@mm@
  \p@tb@x\is@b@x{\sp@c@ no\p@@nt}%
  \ifvoid\ed@b@x\else\sp@c@(\@nb@x\ed@b@x,\space\ifed@ ed.\else eds.\fi)\fi
  \p@tb@x\p@blb@x\c@mm@ \p@tb@x\p@bl@db@x\c@mm@ \p@tb@x\y@@rb@x\c@mm@
  \p@tb@x\p@g@b@x{\c@mm@\ifp@g@ p\p@@nt\else pp\p@@nt\fi}%
  \p@tb@x\inf@b@x\c@mm@\fi
  \else \p@tb@x\p@p@nf@b@x\c@mm@ \p@tb@x\v@lb@x\c@mm@
  \ifvoid\y@@rb@x\else\sp@c@(\@nb@x\y@@rb@x)\fi
  \p@tb@x\is@b@x{\c@mm@ no\p@@nt}\p@tb@x\p@g@b@x\c@mm@ \p@tb@x\inf@b@x\c@mm@\fi
  \else \p@tb@x\b@@kb@x\c@mm@
  \p@tb@x\b@@nf@b@x\c@mm@ \p@tb@x\p@blb@x\c@mm@
  \p@tb@x\p@bl@db@x\c@mm@ \p@tb@x\y@@rb@x\c@mm@
  \ifvoid\p@g@b@x\else\c@mm@\@nb@x\p@g@b@x p\fi \p@tb@x\inf@b@x\c@mm@ \fi
  \else \c@mm@\@nb@x\p@p@rb@x\ic@\p@tb@x\p@p@nf@b@x\c@mm@
  \p@tb@x\v@lb@x\sp@c@ \ifvoid\y@@rb@x\else\sp@c@(\@nb@x\y@@rb@x)\fi
  \p@tb@x\is@b@x{\c@mm@ no\p@@nt}\p@tb@x\p@g@b@x\c@mm@\p@tb@x\inf@b@x\c@mm@\fi
  \else \p@tb@x\p@p@rb@x\c@mm@\ic@\p@tb@x\p@p@nf@b@x\c@mm@
  \c@mm@\@nb@x\inb@@kb@x \p@tb@x\b@@nf@b@x\c@mm@ \p@tb@x\v@lb@x\sp@c@
  \p@tb@x\is@b@x{\sp@c@ no\p@@nt}%
  \ifvoid\ed@b@x\else\sp@c@(\@nb@x\ed@b@x,\space\ifed@ ed.\else eds.\fi)\fi
  \p@tb@x\p@blb@x\c@mm@ \p@tb@x\p@bl@db@x\c@mm@ \p@tb@x\y@@rb@x\c@mm@
  \p@tb@x\p@g@b@x{\c@mm@\ifp@g@ p\p@@nt\else pp\p@@nt\fi}%
  \p@tb@x\inf@b@x\c@mm@\fi
  \else\p@tb@x\p@p@rb@x\c@mm@\ic@\p@tb@x\p@p@nf@b@x\c@mm@\p@tb@x\j@@rb@x\c@mm@
  \p@tb@x\v@lb@x\sp@c@ \ifvoid\y@@rb@x\else\sp@c@(\@nb@x\y@@rb@x)\fi
  \p@tb@x\is@b@x{\c@mm@ no\p@@nt}\p@tb@x\p@g@b@x\c@mm@ \p@tb@x\inf@b@x\c@mm@
 \fi}}
\def\m@r@f#1#2{\endr@f@\ifvoid\p@nct@\else\R@fb@x\hbox{\unhbox\r@f@b@x
 #1\unhbox\p@nct@\penalty-200\enskip#2}\fi\egroup\b@g@nr@f}
\def\endref{\endr@f@\ifvoid\p@nct@\else\R@fb@x\hbox{\unhbox\r@f@b@x.}\fi
 \parindent\r@f@nd
 \r@fit@m{\ifvoid\k@yb@x\else\if@fn@m\r@fn@{\unhbox\k@yb@x}\else
 \if@fl@b\unhbox\k@yb@x\else\r@fk@y{\unhbox\k@yb@x}\fi\fi\fi}\unhbox\r@f@b@x
 \endgraf\egroup\endgroup}
\def\moreref{\m@r@f;\empty}
\def\transl{\m@r@f;{\unskip\space
 {\sl English translation\ic@}:\penalty-66 \space}}
\def\endRefs{\endgraf\goodbreak\endgroup}

\hyphenation{acad-e-my acad-e-mies af-ter-thought anom-aly anom-alies
an-ti-deriv-a-tive an-tin-o-my an-tin-o-mies apoth-e-o-ses
apoth-e-o-sis ap-pen-dix ar-che-typ-al as-sign-a-ble as-sist-ant-ship
as-ymp-tot-ic asyn-chro-nous at-trib-uted at-trib-ut-able bank-rupt
bank-rupt-cy bi-dif-fer-en-tial blue-print busier busiest
cat-a-stroph-ic cat-a-stroph-i-cally con-gress cross-hatched data-base
de-fin-i-tive de-riv-a-tive dis-trib-ute dri-ver dri-vers eco-nom-ics
econ-o-mist elit-ist equi-vari-ant ex-quis-ite ex-tra-or-di-nary
flow-chart for-mi-da-ble forth-right friv-o-lous ge-o-des-ic
ge-o-det-ic geo-met-ric griev-ance griev-ous griev-ous-ly
hexa-dec-i-mal ho-lo-no-my ho-mo-thetic ideals idio-syn-crasy
in-fin-ite-ly in-fin-i-tes-i-mal ir-rev-o-ca-ble key-stroke
lam-en-ta-ble light-weight mal-a-prop-ism man-u-script mar-gin-al
meta-bol-ic me-tab-o-lism meta-lan-guage me-trop-o-lis
met-ro-pol-i-tan mi-nut-est mol-e-cule mono-chrome mono-pole
mo-nop-oly mono-spline mo-not-o-nous mul-ti-fac-eted mul-ti-plic-able
non-euclid-ean non-iso-mor-phic non-smooth par-a-digm par-a-bol-ic
pa-rab-o-loid pa-ram-e-trize para-mount pen-ta-gon phe-nom-e-non
post-script pre-am-ble pro-ce-dur-al pro-hib-i-tive pro-hib-i-tive-ly
pseu-do-dif-fer-en-tial pseu-do-fi-nite pseu-do-nym qua-drat-ic
quad-ra-ture qua-si-smooth qua-si-sta-tion-ary qua-si-tri-an-gu-lar
quin-tes-sence quin-tes-sen-tial re-arrange-ment rec-tan-gle
ret-ri-bu-tion retro-fit retro-fit-ted right-eous right-eous-ness
ro-bot ro-bot-ics sched-ul-ing se-mes-ter semi-def-i-nite
semi-ho-mo-thet-ic set-up se-vere-ly side-step sov-er-eign spe-cious
spher-oid spher-oid-al star-tling star-tling-ly sta-tis-tics
sto-chas-tic straight-est strange-ness strat-a-gem strong-hold
sum-ma-ble symp-to-matic syn-chro-nous topo-graph-i-cal tra-vers-a-ble
tra-ver-sal tra-ver-sals treach-ery turn-around un-at-tached
un-err-ing-ly white-space wide-spread wing-spread wretch-ed
wretch-ed-ly Brown-ian Eng-lish Euler-ian Feb-ru-ary Gauss-ian
Grothen-dieck Hamil-ton-ian Her-mit-ian Jan-u-ary Japan-ese Kor-te-weg
Le-gendre Lip-schitz Lip-schitz-ian Mar-kov-ian Noe-ther-ian
No-vem-ber Rie-mann-ian Schwarz-schild Sep-tem-ber}

\def\leftheadtext#1{} \def\rightheadtext#1{}

\let\nopagenumber\p@gen@false \let\putpagenumber\p@gen@true
\let\pagefirst\nopagenumber \let\pagenext\putpagenumber

\else

\amsppttrue

\let\twelvepoint\relax \let\Twelvepoint\relax \let\putpagenumber\relax
\let\logo@\relax \let\pagefirst\firstpage@true \let\pagenext\firstpage@false
\def\nopagenumber{\let\f@li@ld\folio\def\folio{\global\let\folio\f@li@ld}}

\def\ftext#1{\footnotetext""{\vsk-.8>\nt #1}}

\def\m@R@f@[#1]{\Refs\nofrills{}\m@th\tenpoint
 {
 \s@ct{}{#1}}\wr@@c{\string\Refcd{#1}{\the\pageno}}
 \def\k@yf@##1{\hss[##1]\enspace} \let\keyformat\k@yf@
 \def\widest##1{\s@twd@\refindentwd{\tenpoint\k@yf@{##1}}}
 \let\Key\key \def\refin{\kern\refindentwd}}
\let\info\finalinfo \r@R@fs\Refs
\def\adjust#1{#1} \let\adjustend\relax
\let\adjustnext\adjust 

\fi

\outer\def\myRefs{\myR@fs} \r@st@re\proclaim
\def\bye{\par\vfill\supereject\cl@selbl\cl@secd\b@e} \r@endd@\b@e
\let\Cite\cite \let\Key\key \def\endpro{\par\endproclaim}
\let\d@c@\document \def\document{\d@c@\tenpoint}
\hyphenation{ortho-gon-al}

\newtoks\@@tp@t \@@tp@t\output
\output=\@ft@{\let\{\noexpand\the\@@tp@t}
\let\{\relax

\newif\ifVersion

\def\s@ct#1#2{\ifVersion
 \skip@\lastskip\ifdim\skip@<1.5\bls\vskip-\skip@\p@n@l{-200}\vsk.5>%
 \p@n@l{-200}\vsk.5>\p@n@l{-200}\vsk.5>\p@n@l{-200}\vsk-1.5>\else
 \p@n@l{-200}\fi\ifdim\skip@<.9\bls\vsk.9>\else
 \ifdim\skip@<1.5\bls\vskip\skip@\fi\fi
 \vtop{\twelvepoint\raggedright\bf\vp1\vsk->\vskip.16ex\s@twd@\parindent{#1}%
 \ifdim\parindent>\z@\adv\parindent.5em\fi\hang\textindent{#1}#2\strut}
 \else
 \p@sk@p{-200}{.8\bls}\vtop{\bf\s@twd@\parindent{#1}%
 \ifdim\parindent>\z@\adv\parindent.5em\fi\hang\textindent{#1}#2\strut}\fi
 \nointerlineskip\nobreak\vtop{\strut}\nobreak\vskip-.6\bls\nobreak}

\def\p@n@l#1{\ifnum#1=\z@\else\penalty#1\relax\fi}

\def\s@bs@ct#1#2{\ifVersion
 \skip@\lastskip\ifdim\skip@<1.5\bls\vskip-\skip@\p@n@l{-200}\vsk.5>%
 \p@n@l{-200}\vsk.5>\p@n@l{-200}\vsk.5>\p@n@l{-200}\vsk-1.5>\else
 \p@n@l{-200}\fi\ifdim\skip@<.9\bls\vsk.9>\else
 \ifdim\skip@<1.5\bls\vskip\skip@\fi\fi
 \vtop{\elevenpoint\raggedright\it\vp1\vsk->\vskip.16ex%
 \s@twd@\parindent{#1}\ifdim\parindent>\z@\adv\parindent.5em\fi
 \hang\textindent{#1}#2\strut}
 \else
 \p@sk@p{-200}{.6\bls}\vtop{\it\s@twd@\parindent{#1}%
 \ifdim\parindent>\z@\adv\parindent.5em\fi\hang\textindent{#1}#2\strut}\fi
 \nointerlineskip\nobreak\vtop{\strut}\nobreak\vskip-.8\bls\nobreak}

\def\gadv{\global\adv} \def\gad#1{\gadv#1\@ne} \def\gadneg#1{\gadv#1-\@ne}

\newcount\t@@n \t@@n=10 \newbox\testbox

\newcount\Sno \newcount\Lno \newcount\Fno

\def\pr@cl#1{\r@st@re\pr@c@\pr@c@{#1}\global\let\pr@c@\relax}

\def\tagg#1{\tag"\rlap{\rm(#1)}\kern.01\p@"}
\def\l@L#1{\l@bel{#1}L} \def\l@F#1{\l@bel{#1}F} \def\<#1>{\l@b@l{#1}F}
\def\Tag#1{\tag{\l@F{#1}}} \def\Tagg#1{\tagg{\l@F{#1}}}
\def\Rem{\demo{\sl Remark}} 
\def\Pf#1.{\demo{Proof #1}} \def\epf{\qed\enddemo}
\def\Ap@x{Appendix}
\def\Appendix{\Sno=64 \t@@n\@ne \wr@@c{\string\Appencd}
 \def\sf@rm{\char\the\Sno} \def\sf@rm@{\Ap@x\space\sf@rm} \def\sf@rm@@{\Ap@x}
 \def\s@ct@n##1##2{\s@ct\empty{\setbox\z@\hbox{##1}\ifdim\wd\z@=\z@
 \if##2*\sf@rm@@\else\if##2.\sf@rm@@.\else##2\fi\fi\else
 \if##2*\sf@rm@\else\if##2.\sf@rm@.\else\sf@rm@.\enspace##2\fi\fi\fi}}}
\def\Appcd#1#2#3{\def\Ap@@{\hglue-\l@ftcd\Ap@x}\ifx\@ppl@ne\empty
 \def\l@@b{\@fwd@@{#1}{\space#1}{}}\if*#2\entcd{}{\Ap@@\l@@b}{#3}\else
 \if.#2\entcd{}{\Ap@@\l@@b.}{#3}\else\entcd{}{\Ap@@\l@@b.\enspace#2}{#3}\fi\fi
 \else\def\l@@b{\@fwd@@{#1}{\c@l@b{#1}}{}}\if*#2\entcd{\l@@b}{\Ap@x}{#3}\else
 \if.#2\entcd{\l@@b}{\Ap@x.}{#3}\else\entcd{\l@@b}{#2}{#3}\fi\fi\fi}

\let\s@ct@n\s@ct
\def\s@ct@@[#1]#2{\@ft@\xdef\csname @#1@S@\endcsname{\sf@rm}\wr@@x{}%
 \wr@@x{\string\labeldef{S}\space{\?#1@S?}\space{#1}}%
 {
 \s@ct@n{\sf@rm@}{#2}}\wr@@c{\string\Entcd{\?#1@S?}{#2}{\the\pageno}}}
\def\s@ct@#1{\wr@@x{}{
 \s@ct@n{\sf@rm@}{#1}}\wr@@c{\string\Entcd{\sf@rm}{#1}{\the\pageno}}}
\def\s@ct@e[#1]#2{\@ft@\xdef\csname @#1@S@\endcsname{\sf@rm}\wr@@x{}%
 \wr@@x{\string\labeldef{S}\space{\?#1@S?}\space{#1}}%
 {
 \s@ct@n\empty{#2}}\wr@@c{\string\Entcd{}{#2}{\the\pageno}}}
\def\s@cte#1{\wr@@x{}{
 \s@ct@n\empty{#1}}\wr@@c{\string\Entcd{}{#1}{\the\pageno}}}
\def\theSno#1#2{\dff\?#1@S?{#2}%
 \wr@@x{\string\labeldef{S}\space{#2}\space{#1}}\fi}

\newif\ifd@bn@\d@bn@true
\def\Section{\gad\Sno\ifd@bn@\Fno\z@\Lno\z@\fi\@fn@xt[\s@ct@@\s@ct@}
\def\section{\gad\Sno\ifd@bn@\Fno\z@\Lno\z@\fi\@fn@xt[\s@ct@e\s@cte}
\let\Sect\Section \let\sect\section
\def\subsection{\@fn@xt*\subs@ct@\subs@ct}
\def\subs@ct#1{{
 \s@bs@ct\empty{#1}}\wr@@c{\string\subcd{#1}{\the\pageno}}}
\def\subs@ct@*#1{\vsk->\vsk>{
 \s@bs@ct\empty{#1}}\wr@@c{\string\subcd{#1}{\the\pageno}}}
 \def\Snodef#1{\Sno #1}

\def\l@b@l#1#2{\def\n@@{\csname #2no\endcsname}%
 \if*#1\gad\n@@ \@ft@\xdef\csname @#1@#2@\endcsname{\l@f@rm}\else\def\t@st{#1}%
 \ifx\t@st\empty\gad\n@@ \@ft@\xdef\csname @#1@#2@\endcsname{\l@f@rm}%
 \else\@ft@\ifx\csname @#1@#2@mark\endcsname\relax\gad\n@@
 \@ft@\xdef\csname @#1@#2@\endcsname{\l@f@rm}%
 \@ft@\gdef\csname @#1@#2@mark\endcsname{}%
 \wr@@x{\string\labeldef{#2}\space{\?#1@#2?}\space\ifnum\n@@<10 \space\fi{#1}}%
 \fi\fi\fi}
\def\labeldef#1#2#3{\dff\?#3@#1?{#2}}
\def\Labeldef#1#2#3{\dff\?#3@#1?{#2}\@ft@\gdef\csname @#3@#1@mark\endcsname{}}

\def\l@bel#1#2{\l@b@l{#1}{#2}\?#1@#2?}

\newcount\c@cite
\def\?#1?{\csname @#1@\endcsname}
\def\[{\@fn@xt:\c@t@sect\c@t@}
\def\c@t@#1]{{\c@cite\z@\@fwd@@{\?#1@L?}{\adv\c@cite1}{}%
 \@fwd@@{\?#1@F?}{\adv\c@cite1}{}\@fwd@@{\?#1?}{\adv\c@cite1}{}%
 \relax\ifnum\c@cite=\z@{\bf ???}\wrs@x{No label [#1]}\else
 \ifnum\c@cite=1\let\@@PS\relax\let\@@@\relax\else\let\@@PS\underbar
 \def\@@@{{\rm<}}\fi\@@PS{\?#1?\@@@\?#1@L?\@@@\?#1@F?}\fi}}
\def\(#1){{\rm(\c@t@#1])}}
\def\c@t@s@ct#1{\@fwd@@{\?#1@S?}{\?#1@S?\relax}%
 {{\bf ???}\wrs@x{No section label {#1}}}}
\def\c@t@sect:#1]{\c@t@s@ct{#1}} \let\SNo\c@t@s@ct

\newdimen\l@ftcd \newdimen\r@ghtcd \let\nlc\relax

\def\d@tt@d{\leaders\hbox to 1em{\kern.1em.\hfil}\hfill}
\def\entcd#1#2#3{\item{#1}{#2}\alb\kern.9em\hbox{}\kern-.9em\d@tt@d
 \kern-.36em{#3}\kern-\r@ghtcd\hbox{}\par}
\def\Entcd#1#2#3{\def\l@@b{\@fwd@@{#1}{\c@l@b{#1}}{}}\vsk.2>%
 \entcd{\l@@b}{#2}{#3}}
\def\subcd#1#2{{\adv\leftskip.333em\entcd{}{\it #1}{#2}}}
\def\Refcd#1#2{\def\t@@st{#1}\ifx\t@@st\empty\ifx\r@fl@ne\empty\relax\else
 \R@fcd{\r@fl@ne}{#2}\fi\else\R@fcd{#1}{#2}\fi}
\def\R@fcd#1#2{\sk@@p{.6\bls}\entcd{}{\hglue-\l@ftcd\bf #1}{#2}}
\def\Refline{\def\r@fl@ne} \def\Refempty{\let\r@fl@ne\empty}
\def\Appencd{\par\adv\leftskip-\l@ftcd\adv\rightskip-\r@ghtcd\@ppl@ne
 \adv\leftskip\l@ftcd\adv\rightskip\r@ghtcd\let\Entcd\Appcd}
\def\appline{\def\@ppl@ne} \def\Appempty{\let\@ppl@ne\empty}
\def\Appline#1{\def\@ppl@ne{\s@bs@ct{}{\bf#1}}}
\def\leftcd#1{\adv\leftskip-\l@ftcd\s@twd@\l@ftcd{\c@l@b{#1}\enspace}
 \adv\leftskip\l@ftcd}
\def\rightcd#1{\adv\rightskip-\r@ghtcd\s@twd@\r@ghtcd{#1\enspace}
 \adv\rightskip\r@ghtcd}
\def\C@nt{Contents} \def\Ap@s{Appendices} \def\R@fcs{References}
\def\contents{\@fn@xt*\cont@@\cont@}
\def\cont@{\@fn@xt[\cnt@{\cnt@[\C@nt]}}
\def\cont@@*{\@fn@xt[\cnt@@{\cnt@@[\C@nt]}}
\def\cnt@[#1]{\c@nt@{M}{#1}{44}{\s@bs@ct{}{\bf\Ap@s}}}
\def\cnt@@[#1]{\c@nt@{M}{#1}{44}{}}
\def\endco{\par\penalty-500\vsk>\vskip\z@\endgroup}
\def\readcd{\@np@t{\jobname.cd}}
\def\Cde{\@fn@xt*\Cde@@\Cde@}
\def\Cde@{\@fn@xt[\Cd@{\Cd@[\C@nt]}}
\def\Cde@@*{\@fn@xt[\Cd@@{\Cd@@[\C@nt]}}
\def\Cd@[#1]{\cnt@[#1]\readcd\endco}
\def\Cd@@[#1]{\cnt@@[#1]\readcd\endco}
\def\contlabeldef{\def\c@l@b}

\long\def\c@nt@#1#2#3#4{\s@twd@\l@ftcd{\c@l@b{#1}\enspace}
 \s@twd@\r@ghtcd{#3\enspace}\adv\r@ghtcd1.333em
 \def\@ppl@ne{#4}\def\r@fl@ne{\R@fcs}\s@ct{}{#2}\B@gr@\parindent\z@\let\nlc\nl
 \let\nl\relax\parskip.2\bls\adv\leftskip\l@ftcd\adv\rightskip\r@ghtcd}

\def\writecd{\immediate\openout\@@cd\jobname.cd \def\wr@@c{\write\@@cd}
 \def\cl@secd{\immediate\write\@@cd{\string\endinput}\immediate\closeout\@@cd}
 \def\closecd{\cl@secd\global\let\cl@secd\relax}}
\let\cl@secd\relax \def\wr@@c#1{} \let\closecd\relax

\def\dff{\@ft@\d@f} \def\d@f{\@ft@\def}
\def\edff{\@ft@\ed@f} \def\ed@f{\@ft@\edef}
\def\defi#1#2{\def#1{#2}\wr@@x{\string\def\string#1{#2}}}

\def\qed{\hbox{}\nobreak\hfill\nobreak{\m@th$\,\square$}}
\def\back#1 {\strut\kern-.33em #1\enspace\ignore} 
\def\Text#1{\crcr\noalign{\alb\vsk>\normalbaselines\vsk->\vbox{\nt #1\strut}%
 \nobreak\nointerlineskip\vbox{\strut}\nobreak\vsk->\nobreak}}

\def\hcor#1{\advance\hoffset by #1}
\def\vcor#1{\advance\voffset by #1}
\let\bls\baselineskip \let\ignore\ignorespaces
\ifx\ic@\undefined \let\ic@\/\fi
\def\vsk#1>{\vskip#1\bls} \let\adv\advance
\def\vv#1>{\vadjust{\vsk#1>}\ignore}
\def\vvn#1>{\vadjust{\nobreak\vsk#1>\nobreak}\ignore}
\def\vvv#1>{\vskip\z@\vsk#1>\nt\ignore}
\def\vvgood{\vadjust{\penalty-500}} 
\def\Par{\vsk.5>} \def\setparindent{\edef\Parindent{\the\parindent}}
\def\Type{\vsk.5>\bgroup\parindent\z@\tt\rightskip\z@ plus1em minus1em%
 \spaceskip.3333em \xspaceskip.5em\relax}
\def\endType{\vsk.5>\egroup\nt} 

\let\Hat\widehat \let\Tilde\widetilde \let\dollar\$ \let\ampersand\&
\let\sss\scriptscriptstyle  
\let\vp\vphantom \let\hp\hphantom \let\nt\noindent
\let\cline\centerline \let\lline\leftline \let\rline\rightline
\def\nn#1>{\noalign{\vskip#1\p@@}} \def\NN#1>{\openup#1\p@@}
\def\cnn#1>{\noalign{\vsk#1>}}
\def\Cup{\bigcup\limits} 
\let\Lim\lim \def\lim{\Lim\limits} \let\Sum\sum \def\sum{\Sum\limits}
\def\Plus{\bigoplus\limits} 
\let\Prod\prod \def\prod{\Prod\limits} \let\Int\int \def\int{\Int\limits}

\def\tsum{\mathop{\tsize\Sum}\limits} 
\def\tprod{\mathop{\tsize\Prod}\limits}
\def\&{.\kern.1em} \def\>{{\!\;}} \def\]{{\!\!\;}} \def\){\>\]} \def\}{\]\]}
\def\nl{\leavevmode\hfill\break} \def\~{\leavevmode\@fn@xt~\m@n@s\@md@sh}
\def\m@n@s~{\raise.15ex\mbox{-}} \def\@md@sh{\raise.13ex\hbox{--}}
\let\procent\% \def\%#1{\ifmmode\mathop{#1}\limits\else\procent#1\fi}
\let\@ml@t\" \def\"#1{\ifmmode ^{(#1)}\else\@ml@t#1\fi}
\let\@c@t@\' \def\'#1{\ifmmode _{(#1)}\else\@c@t@#1\fi}
\let\colon\: \def\:{^{\vp|}}

\let\texspace\ \def\ {\ifmmode\alb\fi\texspace}

\let\n@wp@ge\newpage \def\newpage{\endgraf\n@wp@ge}
\let\=\m@th \def\mbox#1{\hbox{\m@th$#1$}}
\def\mtext#1{\text{\m@th$#1$}} \def\^#1{\text{\m@th#1}}
\def\Line#1{\kern-.5\hsize\line{\m@th$\dsize#1$}\kern-.5\hsize}
\def\Lline#1{\kern-.5\hsize\lline{\m@th$\dsize#1$}\kern-.5\hsize}
\def\Cline#1{\kern-.5\hsize\cline{\m@th$\dsize#1$}\kern-.5\hsize}
\def\Rline#1{\kern-.5\hsize\rline{\m@th$\dsize#1$}\kern-.5\hsize}

\def\Ll@p#1{\llap{\m@th$#1$}} \def\Rl@p#1{\rlap{\m@th$#1$}}
 \def\Cl@p#1{\llap{\m@th$#1$\hss}}
\def\Llap#1{\mathchoice{\Ll@p{\dsize#1}}{\Ll@p{\tsize#1}}{\Ll@p{\ssize#1}}%
 {\Ll@p{\sss#1}}}
\def\Clap#1{\mathchoice{\Cl@p{\dsize#1}}{\Cl@p{\tsize#1}}{\Cl@p{\ssize#1}}%
 {\Cl@p{\sss#1}}}
\def\Rlap#1{\mathchoice{\Rl@p{\dsize#1}}{\Rl@p{\tsize#1}}{\Rl@p{\ssize#1}}%
 {\Rl@p{\sss#1}}}
 
\def\LRtph#1#2{\setbox\z@\hbox{#1}\dimen\z@\wd\z@\hbox{\hbox to\dimen\z@{#2}}}
\def\LRph#1#2{\LRtph{\m@th$#1$}{\m@th$#2$}}

\def\Lto#1{\setbox\z@\mbox{\tsize{#1}}%
 \mathrel{\mathop{\hbox to\wd\z@{\rightarrowfill}}\limits#1}}
\def\Lgets#1{\setbox\z@\mbox{\tsize{#1}}%
 \mathrel{\mathop{\hbox to\wd\z@{\leftarrowfill}}\limits#1}}
\def\vpb#1{{\vp{\big(}}^{\]#1}}

\let\alb\allowbreak 
\def\ald{\noalign{\alb}} 

 \let\x\times \let\ox\otimes 
\let\sub\subset  \let\tabs\+
\let\le\leqslant \let\ge\geqslant
\let\der\partial \let\8\infty \let\*\star

\let\map\mapsto  
 
 \def\vert{\ |\ } \def\nin{\not\in}

\let\lb\lbrace \let\rb\rbrace

\def\lsym#1{#1\alb\ldots\relax#1\alb}
\def\lc{\lsym,}  \def\lx{\lsym\x} \def\lox{\lsym\ox}

\def\Re{\mathop{\roman{Re}\>}} \def\Im{\mathop{\roman{Im}\>}}
\def\End{\mathop{\roman{End}\>}} 
\def\im{\mathop{\roman{im}\>}} 
\def\ker{\mathop{\roman{ker}\>}} 
 
\def\tr{\mathop{\roman{tr}}\nolimits}

\def\sgn{\mathop{\roman{sgn}\)}\limits}
 
\def\id{\roman{id}}  
\def\1{^{-1}} \def\_#1{_{\Rlap{#1}}}
\def\vst#1{{\lower1.9\p@@\mbox{\bigr|_{\raise.5\p@@\mbox{\ssize#1}}}}}
\def\vrp#1:#2>{{\vrule height#1 depth#2 width\z@}}
\def\vru#1>{\vrp#1:\z@>} \def\vrd#1>{\vrp\z@:#1>}
\def\qqq{\qquad\quad} 
\def\sscr#1{\raise.3ex\mbox{\sss#1}} \def\@@PS{\bold{OOPS!!!}}
\def\ono{\bigl(1+o(1)\bigr)} 

\def\intcl{\mathop
 {\Rlap{\raise.3ex\mbox{\kern.12em\curvearrowleft}}\int}\limits}
\def\intcr{\mathop
 {\Rlap{\raise.3ex\mbox{\kern.24em\curvearrowright}}\int}\limits}

\def\pms{\raise.25ex\mbox{\ssize\pm}\>}
\def\mps{\raise.25ex\mbox{\ssize\mp}\>}

\let\al\alpha

 \let\Gm\Gamma 
\let\dl\delta \let\Dl\Delta 
 \let\eps\varepsilon \let\epsilon\eps

\let\zt\zeta
 \let\Tht\Theta

\let\si\sigma \let\Si\Sigma
 
\let\pho\phi \let\phi\varphi

\def\C{\Bbb C}
\def\R{\Bbb R}
\def\Z{\Bbb Z}

\def\DD{\Bbb D}
\def\FF{\Bbb F}

\def\Zp{\Z_{\ge 0}} 
\def\Zn{\Z_{\le 0}} 

\def\difl/{differential} \def\dif/{difference}
\def\cf.{cf.\ \ignore} \def\Cf.{Cf.\ \ignore}
\def\egv/{eigenvector} \def\eva/{eigenvalue} \def\eq/{equation}
\def\lhs/{the left hand side} \def\rhs/{the right hand side}
\def\Lhs/{The left hand side} \def\Rhs/{The right hand side}
\def\gby/{generated by} \def\wrt/{with respect to} \def\st/{such that}
\def\resp/{respectively} \def\off/{offdiagonal} \def\wt/{weight}
\def\pol/{polynomial} \def\rat/{rational} \def\tri/{trigonometric}
\def\fn/{function} \def\var/{variable} \def\raf/{\rat/ \fn/}
\def\inv/{invariant} \def\hol/{holomorphic} \def\hof/{\hol/ \fn/}
\def\mer/{meromorphic} \def\mef/{\mer/ \fn/} \def\mult/{multiplicity}
\def\sym/{symmetric} \def\perm/{permutation} \def\fd/{finite-dimensional}
\def\rep/{representation} \def\irr/{irreducible} \def\irrep/{\irr/ \rep/}
\def\hom/{homomorphism} \def\aut/{automorphism} \def\iso/{isomorphism}
\def\lex/{lexicographical} \def\as/{asymptotic} \def\asex/{\as/ expansion}
\def\ndeg/{nondegenerate} \def\neib/{neighbourhood} \def\deq/{\dif/ \eq/}
\def\hw/{highest \wt/} \def\gv/{generating vector} \def\eqv/{equivalent}
\def\msd/{method of steepest descend} \def\pd/{pairwise distinct}
\def\wlg/{without loss of generality} \def\Wlg/{Without loss of generality}
\def\onedim/{one-dimensional} \def\qcl/{quasiclassical}
\def\hgeom/{hyper\-geometric} \def\hint/{\hgeom/ integral}
\def\hwm/{\hw/ module} \def\emod/{evaluation module} \def\Vmod/{Verma module}
\def\symg/{\sym/ group} \def\sol/{solution} \def\eval/{evaluation}
\def\anf/{analytic \fn/} \def\anco/{analytic continuation}
\def\qg/{quantum group} \def\qaff/{quantum affine algebra}

\def\Rm/{\^{$R$-}matrix} \def\Rms/{\^{$R$-}matrices} \def\YB/{Yang-Baxter \eq/}
\def\Ba/{Bethe ansatz} \def\Bv/{Bethe vector} \def\Bae/{\Ba/ \eq/}
\def\KZv/{Knizh\-nik-Zamo\-lod\-chi\-kov} \def\KZvB/{\KZv/-Bernard}
\def\KZ/{{\sl KZ\/}} \def\qKZ/{{\sl qKZ\/}}
\def\KZB/{{\sl KZB\/}} \def\qKZB/{{\sl qKZB\/}}
\def\qKZo/{\qKZ/ operator} \def\qKZc/{\qKZ/ connection}
\def\KZe/{\KZ/ \eq/} \def\qKZe/{\qKZ/ \eq/} \def\qKZBe/{\qKZB/ \eq/}

\def\h@ph{\discretionary{}{}{-}} \def\$#1$-{\,\^{$#1$}\h@ph}

\def\TFT/{Research Insitute for Theoretical Physics}
\def\HY/{University of Helsinki} \def\AoF/{the Academy of Finland}
\def\CNRS/{Supported in part by MAE\~MICECO\~CNRS Fellowship}
\def\LPT/{Laboratoire de Physique Th\'eorique ENSLAPP}
\def\ENSLyon/{\'Ecole Normale Sup\'erieure de Lyon}
\def\LPTaddr/{46, All\'ee d'Italie, 69364 Lyon Cedex 07, France}
\def\enslapp/{URA 14\~36 du CNRS, associ\'ee \`a l'E.N.S.\ de Lyon,
au LAPP d'Annecy et \`a l'Universit\`e de Savoie}
\def\ensemail/{vtarasov\@ enslapp.ens-lyon.fr}
\def\DMS/{Department of Mathematics, Faculty of Science}
\def\DMO/{\DMS/, Osaka University}
\def\DMOaddr/{Toyonaka, Osaka 560, Japan}
\def\dmoemail/{vt\@ math.sci.osaka-u.ac.jp}
\def\SPb/{St\&Peters\-burg}
\def\home/{\SPb/ Branch of Steklov Mathematical Institute}
\def\homeaddr/{Fontanka 27, \SPb/ \,191011, Russia}
\def\homemail/{vt\@ pdmi.ras.ru}
\def\absence/{On leave of absence from \home/}
\def\UNC/{Department of Mathematics, University of North Carolina}
\def\ChH/{Chapel Hill}
\def\UNCaddr/{\ChH/, NC 27599, USA} \def\avemail/{av\@ math.unc.edu}
\def\grant/{NSF grant DMS\~9501290}	
\def\Grant/{Supported in part by \grant/}

\def\Aomoto/{K\&Aomoto}
\def\Dri/{V\]\&G\&Drin\-feld}
\def\Fadd/{L\&D\&Fad\-deev}
\def\Feld/{G\&Felder}
\def\Fre/{I\&B\&Fren\-kel}
\def\Gustaf/{R\&A\&Gustafson}
\def\Kazh/{D\&Kazhdan} \def\Kir/{A\&N\&Kiril\-lov}
\def\Kor/{V\]\&E\&Kore\-pin}
\def\Lusz/{G\&Lusztig}
\def\MN/{M\&Naza\-rov}
\def\Resh/{N\&Reshe\-ti\-khin} \def\Reshy/{N\&\]Yu\&Reshe\-ti\-khin}
\def\SchV/{V\]\&\]V\]\&Schecht\-man} \def\Sch/{V\]\&Schecht\-man}
\def\Skl/{E\&K\&Sklya\-nin}
\def\Smirn/{F\]\&Smirnov} \def\Smirnov/{F\]\&A\&Smirnov}
\def\Takh/{L\&A\&Takh\-tajan}
\def\VT/{V\]\&Ta\-ra\-sov} \def\VoT/{V\]\&O\&Ta\-ra\-sov}
\def\Varch/{A\&\]Var\-chenko} \def\Varn/{A\&N\&\]Var\-chenko}

\def\AMS/{Amer.\ Math.\ Society}
\def\CMP/{Comm.\ Math.\ Phys.{}}
\def\DMJ/{Duke.\ Math.\ J.{}}
\def\Inv/{Invent.\ Math.{}} 
\def\IMRN/{Int.\ Math.\ Res.\ Notices}
\def\JPA/{J.\ Phys.\ A{}}
\def\JSM/{J.\ Soviet\ Math.{}}
\def\LMP/{Lett.\ Math.\ Phys.{}}
\def\LMJ/{Leningrad Math.\ J.{}}
\def\LpMJ/{\SPb/ Math.\ J.{}}
\def\SIAM/{SIAM J.\ Math.\ Anal.{}}
\def\SMNS/{Selecta Math., New Series}
\def\TMP/{Theor.\ Math.\ Phys.{}}
\def\ZNS/{Zap.\ nauch.\ semin. LOMI}

\def\ASMP/{Advanced Series in Math.\ Phys.{}}

\def\AMSa/{AMS \publaddr Providence}
\def\Birk/{Birkh\"auser}
\def\CUP/{Cambridge University Press} \def\CUPa/{\CUP/ \publaddr Cambridge}
\def\Spri/{Springer-Verlag} \def\Spria/{\Spri/ \publaddr Berlin}
\def\WS/{World Scientific} \def\WSa/{\WS/ \publaddr Singapore}

\newbox\lefthbox \newbox\righthbox

\let\sectsep. \let\labelsep. \let\contsep. \let\labelspace\relax
\let\sectpre\relax \let\contpre\relax
\def\sf@rm{\the\Sno} \def\sf@rm@{\sectpre\sf@rm\sectsep}
\def\c@l@b#1{\contpre#1\contsep}
\def\l@f@rm{\ifd@bn@\sf@rm\labelsep\fi\labelspace\the\n@@}

\def\sectformdef{\def\sf@rm}

\let\DoubleNum\d@bn@true \let\SingleNum\d@bn@false

\def\NoNewNum{\let\writeldf\relax\def\l@b@l##1##2{\if*##1%
 \@ft@\xdef\csname @##1@##2@\endcsname{\mbox{*{*}*}}\fi}}
\def\NoNewTime{\def\todaydef##1{\def\today{##1}}
 \def\nowtimedef##1{\def\nowtime{##1}}}
\def\NoInput{\let\Input\input\let\writeldf\relax}
\def\Fixed{\NoNewTime\NoInput}

\tenpoint
\csname beta.def\endcsname
\Fixed

\expandafter\ifx\csname zero.def\endcsname\relax \else\endinput\fi
\expandafter\edef\csname zero.def\endcsname{%
 \catcode`\noexpand\@=\the\catcode`\@\space}

\catcode`\@=11

\newif\ifMPIM

\def\Th#1{\pr@cl{Theorem \l@L{#1}}\ignore}
\def\Lm#1{\pr@cl{Lemma \l@L{#1}}\ignore}
\def\Cr#1{\pr@cl{Corollary \l@L{#1}}\ignore}
\def\Df#1{\pr@cl{Definition \l@L{#1}}\ignore}
\def\Cj#1{\pr@cl{Conjecture \l@L{#1}}\ignore}
\def\Prop#1{\pr@cl{Proposition \l@L{#1}}\ignore}

\def\appendix#1{\Appendix\sect{\Ap@x.\enspace#1}}

\let\h\hbar

\def\Sb{\bold S}
\def\Fc{\Cal F}

\def\Asym{\mathop{\roman{Asym}\)}}
\def\|{{\setbox\z@\mbox{|}\dp\z@.6\dp\z@\lower1.5\p@@\box\z@}}

\def\Cn{\C^{\>n}}  

\def\Cti{\Tilde C}
\def\Kt{\Tilde K}
\def\pht{\Tilde\phi}

\def\Ch{\Hat C}
\def\phh{\hat\phi}
\def\zth{\hat\zt}

\def\Fxl{\Rlap{\,\,\overline{\!\!\phantom\Fc\]}}\Fc^{\ox\ell}}
\def\Fq{\Fc_{\!q}}
\def\Fql{\Rlap{\Fq}{\phantom\Fc}^{\smash{\wedge\ell}}}
\def\Fqb{\Rlap{\,\,\overline{\!\!\phantom\Fc\]}}\Fq}
\def\Fqbw#1{\Rlap{\Fqb}{\phantom\Fc}^{\smash{\wedge#1}}}
\def\Fqbl{\Fqbw\ell}

\def\Ib{\bar I}
\def\Il{I^{\ox\ell}}
\def\Itw{I^{\ox2}}

\def\Sl{\Sb_\ell}

\def\Fp#1{F^{(#1)}}
\def\Xp#1{\Xi^{(#1)}}

\def\Exp(#1){\exp\!\)\bigl(2\pi i\)(#1)/p\)\bigr)}
\def\Gmm(#1){\Gm\bigl((#1)/p\)\bigr)}
\def\'#1{\ifmmode _{\lb#1\rb}\else\@c@t@#1\fi}
\def\oxc{\mathrel{\ox{\vp|}_{\)\C\!}}}

\def\tsun{\tsum_{m=1}^n} 
\def\sun{\sum_{m=1}^n} \def\pron{\prod_{m=1}^n}
\def\pral{\prod_{a=1}^\ell} \def\tpral{\tprod_{a=1}^\ell}
\def\prab{\prod_{1\le a<b\le\ell}}
\def\susi{\sum_{\,\si\in\Sl\!}}

\def\mn{m=1\lc n} \def\aell{a=1\lc\ell}
\def\xin{\xi_1\lc\xi_n} \def\zn{z_1\lc z_n} \def\zmn{z_1\lc z_m\}+p\lc z_n}
\def\xim{\xi_1\lc\xi_{n-1}}
\def\tell{t_1\lc t_\ell} \def\twll{t_2\lc t_\ell} 

\def\Cb#1{\C\)[\)#1\)]} \def\sing{^{\sscr{\italic{sing}}}}
\def\Von{V^{\ox n}} \def\Vl{(\Von)_\ell\:} \def\Vls{\Vns_\ell}
\def\Vns{(\Von)\sing} \def\Vqs{(\Von)^{\sscr{\italic{sing}\ssize,\)q}}}
\def\Vlm{(\Von)_{\ell-1}\:} \def\Vlmm{(\Von)_{\ell-2}\:} \def\Vk{(\Von)_k\:}
\def\Det#1{\mathop{\vp e\smash\det}\limits^{\,\sss\"{#1}\!}}
\def\limq{\lim_{q\to i}} \def\limF{\limq\]\bigl(F(q)\>\Vlm\bigr)}

\def\Imu{0\le\Im\mu<2\pi}
\def\Immu{0<\Im\mu<2\pi}
\def\Rezm{\Re\bigl((z_m\!-z_{m+1})/p\bigr)}

\def\gsl{\frak{sl}_2} \def\Uq{U_q(\gsl)} \def\p-{\^{$\,p\)$}-}

\def\hmap/{\hgeom/ map} \def\hsol/{\hgeom/ \sol/}
\def\fsol/{fundamental \sol/}

\let\goodbm\relax  \let\mmgood\relax

\ifMag\let\goodbm\goodbreak  \let\mmgood\vvgood
 \let\goodbreak\relax  \let\vvgood\relax \fi

\csname zero.def\endcsname

\labeldef{S} {1} {1}

\labeldef{F} {2\labelsep \labelspace 1}  {qKZ}
\labeldef{F} {2\labelsep \labelspace 2}  {Kmz}
\labeldef{F} {2\labelsep \labelspace 3}  {p=2h}

\labeldef{F} {3\labelsep \labelspace 1}  {wM}
\labeldef{F} {3\labelsep \labelspace 2}  {vM}
\labeldef{F} {3\labelsep \labelspace 3}  {hint}
\labeldef{F} {3\labelsep \labelspace 4}  {Df}
\labeldef{L} {3\labelsep \labelspace 2}  {ID}
\labeldef{L} {3\labelsep \labelspace 3}  {ID0}
\labeldef{F} {3\labelsep \labelspace 5}  {PsiW}
\labeldef{L} {3\labelsep \labelspace 4}  {sol}
\labeldef{L} {3\labelsep \labelspace 5}  {sol0}

\labeldef{L} {4\labelsep \labelspace 1}  {compl}
\labeldef{F} {4\labelsep \labelspace 1}  {hmap}
\labeldef{L} {4\labelsep \labelspace 2}  {mu<>0}
\labeldef{L} {4\labelsep \labelspace 3}  {mu=0}
\labeldef{F} {4\labelsep \labelspace 2}  {TXi}
\labeldef{L} {4\labelsep \labelspace 4}  {mu=00}

\labeldef{F} {5\labelsep \labelspace 1}  {WM}
\labeldef{L} {5\labelsep \labelspace 1}  {Wbasis}
\labeldef{L} {5\labelsep \labelspace 2}  {detM}
\labeldef{L} {5\labelsep \labelspace 3}  {detI}
\labeldef{F} {5\labelsep \labelspace 2}  {B}

\labeldef{S} {6} {P0}
\labeldef{L} {6\labelsep \labelspace 1}  {dermu}
\labeldef{F} {6\labelsep \labelspace 1}  {mueq}
\labeldef{F} {6\labelsep \labelspace 2}  {LK}
\labeldef{F} {6\labelsep \labelspace 3}  {AB}
\labeldef{L} {6\labelsep \labelspace 2}  {alpha}
\labeldef{L} {6\labelsep \labelspace 3}  {kerA}

\labeldef{S} {7} {P00}
\labeldef{F} {7\labelsep \labelspace 1}  {phi12}
\labeldef{F} {7\labelsep \labelspace 2}  {BF}
\labeldef{F} {7\labelsep \labelspace 3}  {Xp}
\labeldef{F} {7\labelsep \labelspace 4}  {dim}
\labeldef{F} {7\labelsep \labelspace 5}  {dim2}

\labeldef{S} {8} {Uq}
\labeldef{F} {8\labelsep \labelspace 1}  {JW}
\labeldef{F} {8\labelsep \labelspace 2}  {xiM}
\labeldef{F} {8\labelsep \labelspace 3}  {FFF}

\labeldef{L} {\char 65\labelsep \labelspace 1}	{D1}
\labeldef{F} {\char 65\labelsep \labelspace 1}	{LKt}
\labeldef{F} {\char 65\labelsep \labelspace 2}	{Ktr}

\document
\Versiontrue

\ifMag\hfuzz 15pt\fi

\headline={\hfil\ifnum\pageno=1\hbox{\tenpoint MPI \,98\~~87}\fi}
\center
{\twelvepoint\bf \bls1.2\bls
Completeness of the Hypergeometric Solutions of the \qKZ/ Equation
at Level Zero
\par}
\ifMag\vsk1.2>\else\vsk1.5>\fi
\VT/
\ifMag\vsk1.2>\else\vsk1.4>\fi
{\it \home/
\vsk.1>
\homeaddr/}
\ifMag\vsk1.5>\else\vsk1.62>\fi
{\sl September \,21, 1998}
\ifMag\vsk1.55>\else\vsk1.73>\fi
\endcenter
\ftext{\=\bls12pt 
{\tenpoint\sl E-mail\/{\rm:} 
\homemail/}}

\Abstract
We discuss the \hsol/s of the quantized \KZv/ (\qKZ/) \eq/ at level zero and
show that they give all \sol/s of the \qKZe/. We completely describe linear
relations between the \hsol/s and give an interpretation of the relations
in terms of the \qg/ $\Uq$ at $q=e^{\pi i/2}\!$.
\endAbs
\ifMag\vsk1.25>\else\vsk1.45>\fi

\vsk0>

\Sno 0
\Sect[1]{Introduction}
In this paper we study the \hsol/s of the quantized \KZv/ (\qKZ/) \eq/ at level
zero, \cf. \(qKZ)\,--\,\(p=2h). We are going to address two problems. The first
one is to show that the \hsol/s give all \sol/s of the \qKZe/. The second one
is to describe completely linear relations between the \hsol/s, in other words,
to say in what cases the construction of \hsol/s produces zero \sol/s of
the \qKZe/.
\par
The \qKZe/ \(qKZ) at level zero for $\mu=0$ was introduced in \Cite{S1} as
equations for form factors in the \$\]SU(2)$-invariant Thirring model and a
construction of the \hsol/s was given therein. A simple counting of dimensions
shows that the obtained \hgeom/ \sol/s are not linearly independent. Linear
relations between the \hsol/s were found later in \Cite{S2} as a consequence of
the deformed Riemann bilinear identity, and it was conjectured that there was
no other linear relations between the \hsol/s. This conjecture, in particular,
implied that all \sol/s of the \qKZe/ can be obtained as \hsol/s.
\par
For the \rat/ \$\gsl$-type \qKZe/ at generic level the \hsol/s were contsructed
in \Cite{TV}, \Cite{MV} and it was shown that they exhausted all \sol/s of
the \qKZe/. The general construction was adapted to the case of level zero
for ${\mu=0}$ in \Cite{NPT} and Smirnov's \sol/ were obtained from it,
but completeness of the \hsol/s still remained an open question.
\par
In this paper we consider the \qKZe/ \(qKZ) at level zero both for ${\mu\ne 0}$
and ${\mu=0}$. The completeness of the \hsol/s in the case ${\mu\ne 0}$ is
claimed in Theorem \[mu<>0], the main part of its proof being given by
a formula for a determinant of a matrix formed by the \hint/s, see
Theorem \[detI], which is analogous to Theorem 5.14 in \Cite{TV}.
Theorem \[mu<>0] also implies that in the case $\mu\ne 0$ the \hsol/s
are linearly independent.
\par
To prove the completeness of the \hsol/s in the case ${\mu=0}$, claimed in
Theorem \[mu=0], we employ a new \difl/ \eq/ obeyed by the \hsol/s of
the \qKZe/, see Theorem \[mueq]. Linear relations between the \hsol/s in
the case $\mu=0$ are described in Theorem \[mu=00], its proof being given
in Section \[:P00]. In Section \[:Uq] we interprete of the relations between
the \hsol/s in terms of the \qg/ $\Uq$ at $q=e^{\pi i/2}\!$.
\par
The author thanks A\&Nakayashiki and F\]\&Smirnov for helpful discussions.
The author is grateful to the Max-Planck-Institut f\"ur Mathematik in Bonn,
where the basic part of the paper had been written, for hospitality.

\Sect{The \qKZe/ at level zero}
Let $\h$ be a nonzero complex number. Let ${V=\C\>v_+\oplus\C\>v_-}$
and let ${R(z)\in\End(V\ox V)}$ be the following \Rm/:
$$
R(z)\,=\;{z+\h\)P\over z+\h}\;,
$$
\vsk1.2>\vsk->\nt
where $P\in\End(V\ox V)$ is the permutation operator. Let
$$
\gather
\\
\cnn-1.3>
\si^+\)=\,\left(\){0\ \;1\atop 0\ \;0}\)\right)\,,\qqq
\si^-\)=\,\left(\){0\ \;0\atop 1\ \;0}\)\right)\,,\qqq
\si^3\)=\,\left(\){1\ \)\hp-0\atop 0\ \){-1}}\)\right)
\\
\nn6>
\Text{be the Pauli matrices and}
\nn-8>
H\,=\;{1-\si^3\]\over 2}\;=\,\left(\){0\ \;0\atop 0\ \;1}\)\right)\;.
\endgather
$$
As usual, for any $X\in\End(V)$ we set
$X_i\,=\,\id\lox\%{X}_{\^{$i$-th}}\lox\id$, and similarly,
if $R=\sum r\ox r'\in\End(V\ox V)$, then we set
$$
R_{ij}\,=\>\sum\,\id\lox\%r_{\^{$i$-th}\,}\lox\%{r'}_{\^{$j$-th}\,}\lox\id\,.
$$
\par
Fix a nonzero complex number $p$ called a \em{step} and a complex number $\mu$.
We consider the \qKZe/ for a \$\Von\]$-valued \fn/ $\Psi(\zn)$:
$$
\alignat2
& \Rlap{\Psi(\zmn)\,=\,K_m(\zn)\>\Psi(\zn)\,,\qqq\mn\,,} &&
\Tag{qKZ}
\\
\nn12>
& K_m(\zn)\,&&{}=\, R_{m,m-1}(z_m-z_{m-1}+p)\,\ldots\,R_{m,1}(z_m-z_1+p)\,\x{}
\Tag{Kmz}
\\
\nn6>
&&& {}\>\x\,\exp\](\mu\)H_m)\>R_{m,n}(z_m-z_n)\,\ldots\,
R_{m,m+1}(z_m-z_{m+1})\,.
\endalignat
$$
The number $-2+p/\h$ is called a \em{level} of the \qKZe/. The operators
$K_1\lc K_n$ are called the \qKZ/ operators. They depend on the parameter $\mu$
only via $e^\mu\!$. So \wlg/ we always assume that
$$
\Imu\,.
$$
\par
Let $\,\Si^a=\sun\si^a_m$, $\;a=\pm,3$. The operators $\Si^\pm,\Si^3$ define
\vv.1>
an \$\gsl$-action in $\Von\!$. Denote by $\Vl$ a weight subspace:
$$
\Vl\,=\,\lb\)v\in\Von\vert \Si^3 v\>=\>(n-2\ell)\>v\)\rb\,,
$$
\vvv.3>
which is nontrivial only for $\ell=0,\lc n$,
and by $\Vls\]$ the subspace of singular vectors:
\ifMag\vvn.2>\fi
$$
\Vls\,=\,\lb\)v\in\Vl\vert \Si^+ v\>=\>0\)\rb\,,
$$
\ifMag\vvv.5>\fi
which is nontrivial only for ${2\ell\le n}$. For ${\mu=0}$ the \qKZ/ operators
$K_1\lc K_n$ commute with the \$\gsl$-action in $\Von\!$ while for $\mu\ne 0$
they commute only with the action of $\Si^3$.  Therefore, the space of \sol/s
of the \qKZe/ inherits a weight decomposition of the \$\gsl$-module $\Von\!$
in general, and all the structure of the \$\gsl$-module if $\mu=0$.
\par
In what follows we discuss \sol/s of the \qKZe/ \(qKZ) taking values in
a weight subspace $\Vl$. For $\mu=0$ we write down only \sol/s taking values
in the subspace of singular vectors $\Vls\!$. All other \sol/s can be obtained
from them by the action of $\Si^-\!$.
\Par
In this paper we consider the case of the \qKZe/ at level zero, that is, all
over the paper unless otherwise stated we assume that
\vv-.3>
$$
p\,=\,2\h\,.
\Tag{p=2h}
$$
\vsk-.3>
\vsk0>

\Sect{Hypergeometric \sol/s of the \qKZe/ at level zero}
The \hsol/s of the \qKZe/ at level zero in the case $\mu=0$ are discussed
in detail in \Cite{NPT}. The case ${\mu\ne 0}$ can be considered similarly. For
this reason we give here only the most essential points of the construction.
In general, the \hsol/s of the \rat/ \qKZe/ are obtained in \Cite{TV}.
Here we adapt the general construction from \Cite{TV} to the case $p=2\h$.
\Par
For any \fn/ $f(\tell)$ we set
$$
\Asym f(\tell)\,=\susi\sgn(\si)\>f(t_{\si_1}\lc t_{\si_\ell})\,.
$$
Let $M$ be a subset of $\lb\)1\lc n\)\rb$ \st/ $\#M=\ell$. We write
${M=\lb\)m_1\lc m_\ell\)\rb}$ assuming that $m_1\lsym<m_\ell$.
Introduce a \fn/ $w_M$ as follows:
$$
\gather
\\
\cnn-1.3>
w_M\,=\,\Asym g_M\,,
\Tag{wM}
\\
\nn5>
g_M(\tell)\,=\,\pral\)\Bigl(\,{1\over t_a-z_{m_a}\!}\,\,
\prod_{1\le j<m_a}\!\!{t_a-z_j-\h\over t_a-z_j}\,\Bigr)\!\prab\!(t_a-t_b-\h)\,,
\endgather
$$
and define a vector $v_M\in\Vl$ by the rule:
\vv-.4>
$$
v_M\,=\,v_{\eps_1}\lox v_{\eps_n}\,,\qqq\eps_m\,=\,\biggl\lb
\matrix +\quad\text{if}\ \ m\nin M\\\nn2> -\quad\text{if}\ \ m\in M\endmatrix
\Tag{vM}
$$
\par
Let $\Fqbl$ be a space of \fn/s $F(\tell\);\zn)$ \p-periodic in $\zn$:
$$
\gather
F(\tell\);\zmn)\,=\,F(\tell\);\zn)\,,
\\
\nn6>
\Text{and \st/ the product}
\nn-3>
F(\tell\);\zn)\>\pron\,\pral\bigl(1-\Exp(t_a-z_m)\]\bigr)
\endgather
$$
is an antisymmetric \pol/ in the exponentials
\vv.1>
$e^{2\pi i\)t_1/p}\lc e^{2\pi i\)t_\ell/p}$ of degree less than $n$ in each of
the \var/s. Notice that any $F\in\Fqbl$ is a \p-periodic \fn/ of $\tell$:
$$
F(t_1\lc t_a\}+p\lc t_\ell\);\zn)\,=\,F(\tell\);\zn)\,.
$$
Denote by $\Fql\sub\Fqbl$ a subspace of \fn/s $F$ \st/ the product
$$
F(\tell\);\zn)\>\exp\]\bigl(-2\pi i\}\tsum_{a=1}^\ell\}t_a/p\bigr)\>
\pron\,\pral\bigl(1-\Exp(t_a-z_m)\]\bigr)
$$ is a \pol/ in the exponentials $e^{2\pi i\)t_1/p}\lc e^{2\pi i\)t_\ell/p}$.
\Rem
Here we defined a minimal space $\Fqbl$of periodic functions required to obtain
the \hsol/s of the \qKZe/. In \Cite{NPT} a larger space of periodic functions
is considered. But it is easy to see that it produces the same space of
the \hsol/s of the \qKZe/ \(qKZ) for $\mu=0$.
\enddemo
Let $\pho(t)$ be the \em{phase \fn/}:
\vv-.7>
$$
\pho(t)\,=\,\exp\](\)\mu\)t/p)\>\pron{\Gmm(t-z_m\!-\h)\over\Gmm(t-z_m)}\;.
$$
\Lm{}
For any $\eps>0$ the phase \fn/ $\pho(t)$ has the following \as/s
$$
\pho(t)\,=\,(t/p)^{-n/2}\)\exp\](\)\mu\)t/p)\>\bigl(1+o(1)\bigr)
$$
as \,$t\to\infty$, \ $\eps<|\arg\,(t/p)|<\pi-\eps$.
\endpro
\nt
The statement follows from the Stirling formula.
\Par
Let $C$ be a simple curve separating the sets $\Cup_{m=1}^n(z_m\!+\h+p\)\Zn)$
\vv.1>
and $\Cup_{m=1}^n(z_m\!+p\)\Zp)$, and going from $-ip\)\8$ to $+ip\)\8$.
\vv.1>
More precisely, we assume that the curve $C$ admits a smooth parametrization
$\rho:\R\to\C$ such that $|\)\rho(u)|\to\8$ as $|u|\to\8$ and
$$
{d\over du}\,\rho(u)\,\to\,\al_\pm\,,\qqq\pm\Im(\al_\pm/p)>0\,,\qqq
u\to\pm\8\,,
$$
where $\al_\pm$ are some complex numbers. In principle, the curve $C$ depends
on $\zn$, but for generic $\zn$ one can show that there is the required curve
which does not change under small variations of $\zn$.
\Par
Denote by $\Il(w,W)$ the following integral:
\ifMag\vvn-.5>\fi
$$
\Il(w,W)\,=\;{1\over\ell\)!}\;
\int_{C^\ell}w(\tell)\,W(\tell)\>\tpral\pho(t_a)\,dt_a\,.
\Tag{hint}
$$
\ifMag\vvv-.5>\fi
We are interested in the integrals $\Il(w_M,W)$ where the \fn/ $w_M$,
${\#M=\ell}$, is defined by \(wM) and ${W\in\Fqbl}$, and call them the \hint/s.
If ${\Immu}$, then the integrand of the integral $\Il(w_M,W)$ vanishes
exponentially if $|\)t_1|\lsym+|\)t_\ell|\to\8$ and the integral is clearly
convergent. If $\Im\mu=0$, then the integral $\Il(w_M,W)$ is convergent
provided that either ${2\ell\le n}$ or ${W\in\Fql}\!$. In fact, for arbitrary
$\ell,w_M$ and ${W\in\Fqbl}$ the integral $\Il(w_M,W)$ considered as a \fn/ of
$\mu$ has a finite limit if $\mu$ tends from inside the strip ${\Immu}$ to any
point on the boundary lines $\Im\mu=0$ and $\Im\mu=2\pi$ except $0$ and
$2\pi i$. Further, if necessary, the integrals $\Il(w_M,W)$ should be
understood in the described regularized sense. Notice that the \hint/s do not
depend on a particular choice of the contour $C$.
\Par
The most essential property of the \hint/s is described in Propositions \[ID],
\[ID0].
\Par
Let $D$ be the following operator acting on \fn/s of one \var/:
$$
\gather
\\
\cnn-1.5>
D\]f(t)\,=\,f(t)-e^\mu f(t+p)\>\prod_{j=1}^n {t-z_j-\h\over t-z_j}\;,
\Tag{Df}
\\
\nn-5>
\Text{that is,}
\nn-2>
\pho(t)\,D\]f(t)\,=\,\pho(t)\>f(t)-\pho(t+p)\>f(t+p)\,.
\endgather
$$
\par
Introduce a space $\Fxl$ of \raf/s $f(\tell;\zn)$ \st/ the product
\ifMag\vvn-.5>\fi
$$
f(\tell\);\zn)\>\pron\,\pral\>(t_a-z_m)\pron\,\pral\>(t_a-z_m\}-\h-p)
$$
\ifMag\vvv-.5>\fi
is a \pol/ in $\tell$ of degree less than $2n+\ell$.
\vv.1>
Notice that the space $\Fxl$ contains all the \fn/s $w_M$, $\#M=\ell$.
\vv.1>
The integrals $\Il(f,W)$ for $f\in\Fxl$ and $W\in\Fqbl$ are defined similarly
to the integrals $\Il(w_M,W)$.
\vsk.2>
Denote by $D_a$ the operator $D$ acting \wrt/ the \var/ $t_a$.
\ifMag\vv.1>\fi
We call \fn/s of the form $\sum_{a=1}^\ell\] D_af_a$ \em{total \dif/s},
if $f_a\in\Fxl$ for all $\aell$.
\Prop{ID}
Let ${\Imu}$ and ${\mu\ne 0}$. Then ${I(D_af,W)=0}$, \>$\aell$, \,for any
$f\in\Fxl$ and $\)\botsmash{W\in\Fqbl}\!$.
\endpro
\Prop{ID0}
Let ${\mu=0}$ and either ${W\in\Fqbl}\!$, ${2\ell\le n}$, or ${W\in\Fql}\!$.
Then $I(D_af,W)=0$, \>$\aell$, \,for any $f\in\Fxl$.
\endpro
\nt
The proof of Propositions \[ID], \[ID0] is rather straightforward.
\goodbm
\Par
For any ${W\in\Fqbl}\}$ define a \fn/ $\Psi_W(\zn)$ with values in $\Vl$
as follows:
\ifMag\vvn-.2>\fi
$$
\Psi_W(z_1\lc z_n)\,=\smash{\sum_{\# M=\ell}}\vp\sum\!\Il(w_M,W)\;v_M\,,
\Tag{PsiW}
$$
\ifMag\vvv.5>\fi
where the vectors $v_M$ are given by \(vM).
\Th{sol}
Let ${\Imu}$ and ${\mu\ne 0}$. Then for any ${W\in\Fqbl}\!$ the \fn/ $\Psi_W$
is a \sol/ of the \qKZe/ \(qKZ) taking values in $\Vl$.
\vvgood
\endpro
\Th{sol0}
Let ${\mu=0}$ and either ${W\in\Fqbl}\!$, ${2\ell\le n}$, or ${W\in\Fql}\!$.
Then the \fn/ $\Psi_W$ is a \sol/ of the \qKZe/ \(qKZ) taking values in
$\Vls\!$.
\endpro
\Cr{}
Let ${\mu=0}$, ${2\ell>n}$ and ${W\in\Fql}\!$. Then $\Psi_W=0$.
\endpro
\Pf.
The subspace $\Vls$ is trivial if $2\ell>n$, which implies the claim.
\epf
Theorems \[sol] and \[sol0] follow from general results on formal integral
\rep/s for \sol/s of the \qKZe/, see \Cite{R}\), \Cite{V}, and Propositions
\[ID], \[ID0]. Details of the proof in the case $\mu=0$ are available in
\Cite{NPT}. The proof in the case $\mu\ne 0$ is similar.
\Par
Solutions $\Psi_W$ of the \qKZe/ are called the \em{\hsol/s}.

\Sect{Completeness of the \hsol/s}
Let us formulate accurately the meaning of completeness of the \hsol/s of
the \qKZe/. We will use the following compact notation: $z=(\zn)$.
\par
Consider a vector space $\Cn$ with coordinates $\zn$ and a union of hyperplanes
$\DD\sub\Cn$:
\ifMag\vv-.4>\fi
$$
\DD\,=\,\lb\)z\in\Cn\vert
z_k\]-z_m\!+\h\)\in p\>\Z\,,\quad k,\mn\,,\ \ k\ne m\)\rb\,,
$$
called a \em{discriminant}. Say that $\zn$ are \em{generic} if $(\zn)\nin\DD$.
The complement of the discriminant is a natural definition region for \sol/s
of the \qKZe/ \(qKZ).
\par
Indeed, the \qKZ/ operators $K_1(z)\lc K_n(z)$ are invertible for generic $z$,
since the \Rm/ $R(x)$ is invertible for $x\ne\pm\h$.  Therefore, for generic
$z$ values of a \sol/ $\Psi$ of the \qKZe/ \(qKZ) at points of a lattice
$z+p\)\Z$, are uniquely determined by its value $\Psi(z)$ at the initial point;
the initial value can be an arbitrary vector in $\Von$ because for any $j,\mn$
$$
K_j(\zmn)\>K_m(\zn)\,=\,K_m(z_1\lc z_j\}+p\lc z_n)\>K_j(\zn)\,.
\ifMag\vp{\Big|}\fi
$$
\par
The contour $C$ used in the definition of the \hint/ \(hint) exists for generic
$z$. Therefore, if a periodic \fn/ ${W\in\Fqbl}$ is regular in a vicinity of
a generic point $z$, so is the corresponding \hsol/ $\Psi_W$ of the \qKZe/.
Hence, we can formulate completeness of the \hsol/s of the \qKZe/ in
the following way.
\Th{compl}
Let $\zn$ be generic. If $\,\Imu$ and $\mu\ne 0$, then
$$
\gather
\lb\)\Psi_W(\zn)\vert W\in\Fqbl\)\rb\,=\,\Vl\,.
\\
\ifMag\nn3>\fi
\Text{If $\,\mu=0\)$ and $\,2\ell\le n$, then}
\ifMag\nn3>\fi
\lb\)\Psi_W(\zn)\vert W\in\Fqbl\)\rb\,=\,\Vls.
\endgather
$$
\endpro
Given a point ${z\in\Cn}$ consider a space $\Fqbl(z)$ of functions in $\tell$
which have the form $W(\cdot\>;z)$ for some ${W\in\Fqbl}$. In other words,
$\Fqbl(z)$ is a space of \fn/s $F(\tell)$ \st/ the product
$$
F(\tell)\>\pron\,\pral\bigl(1-\Exp(t_a-z_m)\]\bigr)
$$
is an antisymmetric \pol/ in the exponentials
$e^{2\pi i\)t_1/p}\lc e^{2\pi i\)t_\ell/p}$ of degree less than $n$
\vv.1>
in each of the \var/s. The space $\Fqbl(z)$ is \fd/ and
$\dsize\dim\)\Fqbl(z)\)=\)\dim\]\Vl\)={n\choose \ell}$.
\par
Define the \em{\hmap/} $\Ib_{z,\ell}:\)\Fqbl(z)\to\Vl$ by the rule:
\ifMag\vv.2>\fi
$$
\Ib_{z,\ell}\)W\,=\sum_{\# M=\ell}\!\Il(w_M,W)\;v_M\,,
\Tag{hmap}
$$
\ifMag\else\vvv-.2>\fi
so that $\dsize\Psi_W(z)\)=\)\Ib_{z,\ell}\)W(\cdot\>;z)$, \cf. \(PsiW).
Then we write down Theorem \[compl] in terms of the \hmap/.
\Th{mu<>0}
Let $\mu\ne 0$. Then for generic $\zn$ the \hmap/
$\Ib_{z,\ell}:\)\Fqbl(z)\to\Vl$ is bijective.
\endpro
\nt
The proof is given in the next section.
\Th{mu=0}
Let $\mu=0$. Then for generic $\zn$ we have $\>\im\Ib_{z,\ell}=\Vls\!$.
\endpro
\nt
The proof is given in Section \[:P0].
\Par
For ${\mu\ne0}$ the \hmap/ $\Ib_{z,\ell}$ has a trivial kernel. So any nonzero
periodic \fn/ ${W\in\Fqbl}\!$ produces a nonzero \sol/ $\Psi_W$, that is,
the space of \sol/s of the \qKZe/ \(qKZ) with values in $\Vl$ is isomorphic to
the space $\Fqbl\!$ of periodic \fn/s. On the contrary, for $\mu=0$ the \hmap/
$\Ib_{z,\ell}$ has a nontrivial kernel, which means that for certain periodic
\fn/s the corresponding \hsol/s vanish, and we have a problem of describing
those periodic \fn/s which produces zero \sol/s of the \qKZe/.
\Par
Introduce the following \fn/s:
$$
\gather
\\
\cnn-1.5>
\Xp1(t)\,=\,\Tht(t)-1\,,\qqq
\Tht(t)\,=\,\pron\,{\Exp(t-z_m)+1\over\Exp(t-z_m)-1}\;,
\Tag{TXi}
\\
\ifMag\nn12>\else\nn10>\fi
\Xp2(t_1,t_2)\,=\,\bigl(\Tht(t_1)\>\Tht(t_2)-1\bigr)\,
{\Exp(t_1-t_2)-1\over\Exp(t_1-t_2)+1}\,+\>\Tht(t_1)-\Tht(t_2)\,.
\endgather
$$
It is easy to see that $\Xp1\in\Fqbw1$ and $\Xp2\in\Fqbw2\!$.
Define linear maps
\ifMag
$$
\gather
X\"1_{z,\ell}\):\>\Fqbw{(\ell-1)}(z)\,\to\,\Fqbl(z)\,,
\\
\cnn2.15>
X\"2_{z,\ell}\):\>\Fqbw{(\ell-2)}(z)\,\to\,\Fqbl(z)\,,
\\
\cnn-3.15>
{\align
& X\"1_{z,\ell}\> F(\tell)\,=\,\Asym\bigl(\)\Xp1(t_1)\>F(\twll)\bigr)\,,
\\
\cnn2.25>
& \Rlap{X\"2_{z,\ell}\>F(\tell)\,=\,
\Asym\bigl(\)\Xp2(t_1,t_2)\>F(t_3\lc t_\ell)\bigr)\,.}
\endalign}
\endgather
$$
\else
$$
\align
X\"1_{z,\ell}\):\>\Fqbw{(\ell-1)}(z)\,\to\,\Fqbl(z)\,,\qqq &
X\"1_{z,\ell}\> F(\tell)\,=\,\Asym\bigl(\)\Xp1(t_1)\>F(\twll)\bigr)\,,
\\
\nn6>
X\"2_{z,\ell}\):\>\Fqbw{(\ell-2)}(z)\,\to\,\Fqbl(z)\,,\qqq &
\Rlap{X\"2_{z,\ell}\>F(\tell)\,=\,
\Asym\bigl(\)\Xp2(t_1,t_2)\>F(t_3\lc t_\ell)\bigr)\,.}
\endalign
$$
\fi
They obviously induce linear maps $X\"1_\ell\!:\>\Fqbw{(\ell-1)}\to\,\Fqbl\]$
and $X\"2_\ell\!:\>\Fqbw{(\ell-2)}\to\,\Fqbl\}$.
\Th{mu=00}
Let $\mu=0$ and $2\ell\le n$. Then for generic $\zn$ we have
${\ker\Ib_{z,\ell}\)={}}\alb{\)\im X\"1_{z,\ell}+{}}\alb\)\im X\"2_{z,\ell}$.
Here we assume by convention that $\>\im X\"a_{z,\ell}\)=\)0\)$ for $\ell<a$.
\endpro
\nt
The proof is given in Section \[:P00].
\Par
The last theorem means that if $\mu=0$ and a periodic \fn/ $W$ has the form
$W=X\"1_\ell W_1+X\"2_\ell W_2$ where $W_1,W_2$ are some periodic \fn/s in less
numbers of \var/s, then $\Psi_W=0$, and vice versa, that is, the space of
\sol/s of the \qKZe/ \(qKZ) at $\mu=0$ with values in $\Vls$ is isomorphic to
the quotient space $\Fqbl\]\big/\bigl(\im X\"1_\ell\!\}+\im X\"2_\ell\bigr)$.
An interesting question is to see whether there is a distinguished subspace
in $\Fqbl\!$ complementary to $\im X\"1_\ell\!\}+\im X\"2_\ell\!$.
\Rem
The relation $\Ib_{z,\ell}\>X\"2_{z,\ell}\>W=0$ was originally discovered by
F\]\&Smirnov in \Cite{S2} as a consequence of the deformed Riemann bilinear
identity.  In Section \[:Uq] we give another interpretation of this relation
in terms of the representation theory of the quantum group $\Uq$ at
$q=e^{\pi i/2}\!$.
\enddemo

\Sect{Proof of Theorem \[mu<>0]}
Recall that we assume $p=2\h$, so that $\zn$ are generic if and only if
$$
\prod_{1\le k<m\le n}(e^{2\pi iz_k/p}+e^{2\pi iz_m/p})\,\ne\,0.
$$
\par
Let $M$ be a subset of $\lb\)1\lc n\)\rb$, ${M=\lb\)m_1\lsym<m_\ell\)\rb}$.
Introduce a \fn/ $W_M(\tell;\zn)$ as follows:
$$
\gather
\\
\ifMag\cnn->\else\cnn-1.3>\fi
W_M\,=\,\Asym G_M\,,
\Tag{WM}
\\
\nn5>
G_M(\tell)\,=\,\pral\)\Bigl(\,\){1\over\Exp(t_a-z_{m_a})-1\!}\,\,
\prod_{1\le j<m_a}\!\!{\Exp(t_a-z_j)+1\over\Exp(t_a-z_j)-1}\,\Bigr)\,,
\endgather
$$
\cf. \(wM). It is clear that $W_M\in\Fqbl$.
\Prop{Wbasis}
For generic $\zn$ \fn/s $W_M(\cdot\>,z)$, $\#M=\ell$, form a basis in
the space $\Fqbl(z)$.
\endpro
\Pf.
First observe that it suffices to prove the statement for $\ell=1$.
For $\ell=1$ the statement is implied by Lemma \[detM] after a substitution:
$u=e^{2\pi i\)t/p}$, $x_m=y_m=-\)e^{2\pi iz_m/p}$, $\mn$.
\mmgood
\epf
\Lm{detM}
Define a matrix $M$ as follows: $\sum_{k=1}^n M_{jk}\>u^{k-1}\>=\!
\prod_{1\le i<j}\}(u-x_i)\prod_{j<i\le n}\}(u+y_i)$. Then
$$
\det M\,=\prod_{1\le i<j\le n}(x_i\]+y_j)\,.
$$
\endpro
\nt
The proof is similar to the proof of the Vandermonde determinant.
\Th{detI}
Let $\Imu$ and $\zn$ be generic. Let either $\mu\ne 0$ or $2\ell\le n$. Then
\ifMag
$$
\alignat2
\\
\cnn-.9>
\det\] &\) \bigl[\Il(w_M,W_N)\bigr]_{\vtop{\bls.7\bls
\mbox{\ssize M,N\sub\lb\)1\lc n\)\rb}\mbox{\ssize\#M=\#N=\ell}}}\,={} &&
\\
\nn4>
&\!\]{}=\Bigl((ip)^{n(n-1)/2}\>\bigl(i\>\Gm(-1/2)\bigr)^n\,(e^\mu\}-1)\vpb{n/2}
\>\exp\}\bigl(\>\mu\!\tsun\! z_m/p\)\bigr)\Bigr)^{\tsize{\)n-1\)\choose\ell-1}}
&& {}\!\!\x{}
\\
\nn4>
&& \Llap{{}\x\,\prod_{1\le k<m\le n}\!\!(z_k-z_m\!-\h)
{\vp{\Big|}}^{\!\!{-}\tsize{\)n-2\)\choose\ell-1}}} &.
\\
\cnn-.4>
\endalignat
$$
\else
$$
\align
\\
\cnn-.9>
\det\] &\) \bigl[\Il(w_M,W_N)\bigr]_{\vtop{\bls.7\bls
\mbox{\ssize M,N\sub\lb\)1\lc n\)\rb}\mbox{\ssize\#M=\#N=\ell}}}\,={}
\\
\nn4>
&\!\]{}=\Bigl((ip)^{n(n-1)/2}\>\bigl(i\>\Gm(-1/2)\bigr)^n\,(e^\mu\}-1)\vpb{n/2}
\>\exp\}\bigl(\>\mu\!\tsun\! z_m/p\)\bigr)\Bigr)^{\tsize{\)n-1\)\choose\ell-1}}
\!\!\prod_{1\le k<m\le n}\!\!(z_k-z_m\!-\h)
{\vp{\Big|}}^{\!\!{-}\tsize{\)n-2\)\choose\ell-1}}.\kern-1em
\endalign
$$
\fi
Here $\)0\le\arg\)(e^\mu\}-1)<2\pi$.
\endpro
\Pf {\rm(}idea\/{\rm)}.
The proof of the theorem is in common with the proof of Theorem 5.14
in \Cite{TV}. Let us indicate here only basic points.
\par
Denote by $D(z)$ the determinant in question and by $E(z)$ \rhs/ of the formula
to be proved. Clearly, it is enough to establish the equality $D(z)=E(z)$ under
the assumption $\Immu$, since both $D(z)$ and $E(z)$ are continuous up to
the boundary line $\Im\mu=0$.
\vsk.3>
The determinant $D(z)$ obeys a system of \deq/s following from the \qKZe/:
$$
\gather
D(\zmn)\,=\,\Det\ell K_m(\zn)\>D(\zn)\,,\qqq\mn\,,\kern-3em
\\
\nn6>
\Text{where}
\nn-6>
\Det\ell K_m\,=\,e{\vp{\Big|}}^{\mu\){\tsize{\)n-1\)\choose\ell-1}}}\>\Bigl(\)
\prod_{1\le j<m}\>{z_m\!-z_j-\h+p\over z_m\!-z_j+\h+p}\;
\prod_{m<j\le n}\>{z_m\!-z_j-\h\over z_m\!-z_j+\h}\;
\Bigr)^{\tsize{\)n-2\)\choose\ell-1}}
\endgather
$$
is a determinant of the restriction of the operator $K_m$ to the weight
subspace $\Vl$. This shows that the ratio $D(z)/E(z)$ is a \p-periodic \fn/
of $\zn$ (quasiconstant).
\vsk.1>
To find this quasiconstant we look at the behaviour of the \hint/s
$\Il(w_M,W_N)$ as $\Rezm$ tends to $+\8$ for all $\mn-1$ and compute the \as/s
of the determinant $D(z)$ in this limit. As a result, we obtain that the ratio
$D(z)/E(z)$ is a quasiconstant which tends to $1$ as $\Rezm\to+\8$ for all
$\mn-1$. Such a quasiconstant necesserily equals $1$, which yields the theorem.
\vsk.33>
The method to determine the required \as/s of the \hint/s as $\Rezm\to+\8$ is
developed in \Cite{TV}\). We describe it briefly for the case in question.
\par
Take a simple curve $\Cti$ going from $-i\8$ to $+i\8$, separating the sets
$\Zp$ and $1/2+\Zn$ and \st/ $|\}\Re u\)|<1$ for any $u\in\Cti$. For any
$x\in\C$ define $\Ch(x)=x+p\>\Cti$.
\vsk.3>
It is clear that $\Il(w_M,W_N)=\ell\)!\,\Il(w_M,G_N)$. The integrand of
the integral $\Il(w_M,G_N)$ has poles at hyperplanes $t_a=z_j+p\)s$, $s\in\Zp$,
\vv.1>
only for $j\le n_a$, and at hyperplanes $t_a=z_j+\h+p\)s$, ${s\in\Zn}$, only
for ${j\ge n_a}$. Therefore, if ${\Re\bigl((z_m\!-z_{m+1})/p\bigr)}$ is large
positive for any ${\mn-1}$, then the contour $C^\ell$ in the definition of
the integral $\Il(w_M,G_N)$ can be replaced by the contour
$\Ch^\ell(z)=\Ch(z_{n_1})\lx\Ch(z_{n_\ell})$ without changing the integral,
so we get
$$
\Il(w_M,W_N)\,=\!
\int_{\Ch^\ell(z)}\!\!w_M(\tell)\,G_N(\tell)\>\tpral\pho(t_a)\,dt_a\,,
$$
\vvv-.3>
\cf. \(hint). Substituting into the last integral the explicit expressions for
$w_M$, $G_N$ and $\pho$, we observe that in the limit $\Rezm\to+\8$, $\mn-1$,
most of the terms of the integrand can be replaced by their \as/s$\]^*\!\}$.%
\ftext{\fitem{$\!^*\]$}
It is important at this moment that we assume ${\Immu}$, so the integrand
decays exponentially as $|\)t_1|\lsym+|\)t_\ell|\to\8$.}
As a result, the integral splits into a sum of products of \onedim/ integrals.
The \onedim/ integrals can be calculated explicitly using the formulae
$$
\int_{\!\]\Cti\,\)} e^{u\>(\mu-\pi i)}\)\Gm(u-1/2)\>\Gm(k-u)\,du\,=\,
2\pi\,\Gm(k-1/2)\;e^{k\mu}\,(e^\mu\}-1)\vpb{1/2-k},
\ifMag\else\qqq k\in\Zp\,.\kern-3em\fi
\Tag{B}
$$
\vvv-.4>
\ifMag$k\in\Zp$, \fi where $\)0\le\arg\)(e^\mu\}-1)<2\pi$.
The formulae follow from the formula for the Barnes integral \Cite{WW}.
\par
At last, the leading term of the \as/s of the integral $\Il(w_M,W_N)$
takes the form
$$
\exp\}\bigl(\>\mu\]\tsum_{a=1}^\ell\]z_{n_a}/p\bigr)\!
\prod_{1\le j<k\le n}\!\!(z_j-z_k)\vpb{d_{jkMN}\:}
\bigl(c_{MN}\:\>(e^\mu\}-1)\vpb{r_{MN}\:}\]+\>o(1)\bigr)
$$
\vvv-.4>
for a certain coefficient $c_{MN}\:$ and exponents $d_{ijMN}\:$, $r_{MN}\:$,
whereof one can see that $\Il(w_M,W_N)=o\bigl(\Il(w_N,W_N)\bigr)$ for $M\ne N$.
Hence, for the determinant $D(z)$ we have
$$
D(z)\,=\prod_{\#M=\ell}\Il(w_M,W_M)\;\ono
$$
\vvv-.4>
as $\Rezm\to+\8$, $\mn-1$. The remaining computation
\vv.1>
which shows that $D(z)/E(z)\to 1$ in this limit is rather straightforward.
\epf
\Rem
Below we describe a \difl/ \eq/ \wrt/ $\mu$ obeyed by the \hsol/s of the
\qKZe/, see \(mueq). It implies a \difl/ \eq/ \wrt/ $\mu$ for the determinant
$\det\bigl[\Il(w_M,W_N)\bigr]$, which allows to give another proof of
Theorem \[detI] using the \difl/ \eq/ and \as/s of the \hint/s $\Il(w_M,W_N)$
as $\Re\mu\to\pm\8$.
\goodbm
\vsk.2>
To illustrate the idea let us sketch a proof of formula \(B).
Both sides of the formula satisfy a \difl/ \eq/
$$
2\>(e^\mu\}-1)\,\)\der_\mu\)f(\mu)\,=\,(e^\mu\}-2k)\>f(\mu)\,,
$$
\vvv.3>
hence they are proportional. The leading term of \as/s of the integral in \lhs/
as $\Re\mu\to+\8$ is determined by the rightmost pole of the integrand lying to
the left of the integration contour; the pole is located at $u=1/2$.
The \as/s of the integral equals $2\pi\)\Gm(k-1/2)\,e^{\mu/2}\>\ono$ and
clearly coincides with the \as/s of \rhs/. Formula \(B) is proved.

\enddemo
Theorem \[mu<>0] follows from definition \(hmap)  of the \hmap/,
Proposition \[Wbasis] and Theorem \[detI].

\Sect[P0]{Proof of Theorem \[mu=0]}
Our proof of Theorem \[mu=0] is based on an additional \difl/ \eq/ which holds
for \hsol/s of the \qKZe/.
\Th{dermu}
Let $\Psi$ be a \hsol/ of the \qKZe/ \(qKZ), see \(PsiW). Then
$$
p\,\der_\mu\Psi\,=\,\Bigl(\)\tsun\!z_m\>H_m\,+\;
{\h\over e^\mu\}-1}\,\bigl(\)e^\mu\!\]\tsun\! H_m +\!\tsum_{1\le k<m\le n}\!
(e^\mu\>\si^-_k\si^+_m+\si^+_k\si^-_m\))\bigr)\Bigr)\>\Psi\,.
\Tag{mueq}
$$
\endpro
\nt
The proof of the theorem is given in Appendix. The theorem is a particular
specialization of a more general result about \hsol/s of the \qKZe/s and
the \KZv/ differential \eq/s \Cite{FTV}.
\Rem
One can check directly that the \eq/ \(mueq) is compatible with
the \qKZe/ \(qKZ):
\ifMag
$$
\align
L(\zmn)\>K_m(\zn)\,&{}=\,p\,\der_\mu K_m(\zn)\,+{}
\Tagg{LK}
\\
\nn6>
&\>{}+\,K_m(\zn)\>L(\zn)\,,\kern-1em
\endalign
$$
\else
$$
L(\zmn)\>K_m(\zn)\,=\,p\,\der_\mu K_m(\zn)\)+K_m(\zn)\>L(\zn)\,,
\Tag{LK}
$$
\fi
$\mn$, where $L$ is the operator in \rhs/ of \eq/ \(mueq), see Appendix.
\vvgood
\enddemo
Let $U$ be a \fd/ vector space and $A(\mu)$ be an \$\End(U)$-valued \fn/ \hol/
in a \neib/ of $\mu=0$. Consider a \difl/ \eq/:
\ifMag\vvn.2>\fi
$$
\mu\>\der_\mu\Psi\,=\,A(\mu)\>\Psi\,.
\Tag{AB}
$$
\ifMag\vvv.2>\fi
Set ${A_0=A(0)}$. Assume that the operator $A_0$ is diagonalizable,
and let ${\psi\in U}$ be an \egv/ of $A_0$: $A_0\)\psi=\al\)\psi$.
\Prop{alpha}
There is a \sol/ $\Psi$ of the \eq/ \(AB) which have the form
$$
\Psi(\mu)\,=\,\mu^\al\)\bigl(\>\psi\>+ \tsum_{k=1}^\8\>\tsum_{j=0}^{i_k}
\mu^k\)(\log\mu)^j\> \psi_{jk}\bigr)
$$
where $0=i_0\le i_1\le i_2\le\ldots$ are integers \st/ ${i_k+1\ge i_{k+1}}$
for any ${k\ge 0}$ and $i_k=i_{k+1}$ for large $k$, and $\psi_{jk}\in U$.
The series is convergent for $\mu$ in a punctured \neib/ of $\mu=0$.
\mmgood
\endpro
\nt
The proof is an exercise on the analytic theory of \difl/ \eq/s;
we leave it to a reader.
\Cr{kerA}
Assume that all nonzero \eva/s of $A_0$ have positive real parts. Let
$\Psi(\mu)\in\End(U)$ be a \fsol/ of the \eq/ \(AB), $\det\Psi(\mu)\ne 0$,
in a punctured \neib/ of $\mu=0$. Then there is a limit
$\Psi_0=\lim_{\mu\to 0}\Psi(\mu)$ and $\im\Psi_0=\ker A_0$.
\endpro
\Pf.
According to Proposition \[alpha] we can start from an eigenbasis of $A_0$ in
$U$ and lift it to a \fsol/ of the \eq/ \(AB). For this \fsol/ the claim of
the corollary clearly holds. Since two \fsol/s of the same \difl/ \eq/ differ
only by multiplication from the right by a \ndeg/ constant matrix, the
corollary follows.
\epf
For the \eq/ \(mueq) taking into account that $p=2\h$ and
$H_m=\si^-_m\)\si^+_m$ we have
$$
A_0\,=\,{1\over 2}\,\Si^-\>\Si^+\)=\,
{1\over 8}\,\bigl(\)\Si^3\)(\)\Si^3\!+2\))+4\>\Si^-\>\Si^+\)\bigr)-
{1\over 8}\,\Si^3\)(\)\Si^3\!+2\))
$$
where $\,\topsmash{\Si^a=\sun\si^a_m}$, $\;a=\pm,3$. Since the first term of
the last expression is the Casimir operator for $\gsl$, the operator $A_0$
restricted to the weight subspace $\Vl$ has the following \eva/s:
$$
{(n-2k)\>(n-2k+2)-(n-2\ell)\>(n-2\ell+2)\over 8}\;=\;
{(\ell-k)\>(n-k-\ell+1)\over 2}\;,
$$
$k=0\lc\min\)(\ell,n-\ell\))$, which are nonnegative half-integers.
Besides, $\ker A_0\cap\Vl=\Vls\!$.
\Par
Consider the \hmap/ $\Ib_{z,\ell}\)(\mu)$, \cf. \(hmap), where we write its
dependence of $\mu$ explicitly. Theorem \[dermu] means that
$\Ib_{z,\ell}\)(\mu)$ obeys the \eq/ \(mueq), and according to Theorem \[mu<>0]
\vv-.08>
it is a \fsol/ of the \eq/ \(mueq) with values in $\Vl$. Hence, by Corollary
\[kerA] there is a limit
\vv-.07>
$\lim_{\mu\to 0}\Ib_{z,\ell}\)(\mu)$ and
$\im\}\bigl(\lim_{\mu\to 0}\Ib_{z,\ell}\)(\mu)\bigr)=\Vls\!$.
If $2\ell\le n$, then $\lim_{\mu\to 0}\Ib_{z,\ell}\)(\mu)=\Ib_{z,\ell}\)(0)$,
which completes the proof of Theorem \[mu=0].
\Rem
The fact that $\)\im\Ib_{z,\ell}\)(0)\sub\Vls$ is well known,
see \Cite{NPT} for details of its proof.
\enddemo

\Sect[P00]{Proof of Theorem \[mu=00]}
It is known that $\>\im X\"1_{z,\ell}\sub\ker\Ib_{z,\ell}$,
see Lemma 5.3 in \Cite{NPT} for details.
\vv-.15>
The proof of the second inclusion \>$\im X\"2_{z,\ell}\sub\ker\Ib_{z,\ell}$
is similar, though slightly more complicated, see Appendix.
Because of Theorem \[mu=0] what remains to show is
$$
\dim\]\bigl(\]\im X\"1_{z,\ell}+\]\im X\"2_{z,\ell}\>\bigr)\,=\,
\dim\]\Vl-\]\dim\]\Vls=\,\dim\)\Si^-\Vlm.
$$
The second equality is due to the common knowledge:
\ifMag
$$
\Vl\,=\,\Vls\oplus\>\Si^-\Vlm.
$$
\par
\else
$\Vl=\Vls\oplus\>\Si^-\Vlm\!$.
\Par
\fi
Let $\xin$ be the Grassman \var/s: $\xi_k\>\xi_m=-\)\xi_m\)\xi_k\>$ for any
$k,\mn$. Consider the group algebra $\Cb\xin$. Set
$\deg\)\xi_1\}\lsym=\deg\)\xi_n\}=1$, and denote by $\Cb\xin_k\:$ the subspace
of elements of degree $k$. By convention, set $\Cb\xin_k\:=0$ for $k<0$.
Introduce elements $\phi\"1\},\)\phi\"2\in\Cb\xin$:
\ifMag\vvn-.3>\else\vvn-.7>\fi
$$
\phi\"1=\tsun\xi_m\,,\qquad\qquad\phi\"2=\tsum_{1\le k<m\le n}\!\xi_k\>\xi_m\,.
\Tag{phi12}
$$
\vvv-.3>
Consider \iso/s of vector spaces:
$$
\Fqbw k(z)\,\to\,\Cb\xin_k\:\,,\qquad\qquad
W_M(\cdot\>,z)\,\map\,\xi_{m_1}\ldots\xi_{m_k}\,,
\Tag{BF}
$$
where ${M=\lb\)m_1\lsym<m_k\)\rb}$. Under these \iso/s the operator
$X\"a_{z,\ell}$ translates modulo proportionality into multiplication
by $\phi\"a\!$ acting on $\Cb\xin_{(\ell-a)}\:$, because the \fn/s
$\Xp1,\,\Xp2\!$, \cf. \(TXi), can be written as follows:
\vvn-.2>
$$
\Xp1=\,2\tsun\}W\'m\,,\qquad\qquad\Xp2=\,4\tsum_{1\le k<m\le n}\!\!W\'{k,m}\,,
\Tag{Xp}
$$
the proof beeing given at the end of the section. Here the \fn/s $W\'m$ and
$W\'{k,m}$ are given by \(WM) for $M=\lb\)m\rb$ and $M=\lb\)k,m\rb$, \resp/.
\Par
Therefore, the problem is to show that
$$
\dim\]\bigl(\)\phi\"1\)\Cb\xin_{(\ell-1)}\:+\>
\phi\"2\)\Cb\xin_{(\ell-2)}\:\)\bigr)\,=\,\dim\)\Si^-\Vlm.
\Tag{dim}
$$
The cases $\ell=0,1$ are trivial. The case $\ell=n=2$ is simple. From now on
we assume that $\ell\ge 2$ and $n\ge 3$ (\)recall that $\ell\le n$)\).
It is easy to see that
$$
\gather
\\
\cnn-1.2>
\Cb\xin_k\:\,=\,\Cb\xim_k\:\>\oplus\>\phi\"1\Cb\xim_{k-1}\:
\\
\nn-2>
\Text{and}
\nn-6>
\phi\"1\Cb\xim_{k-1}\:\,=\,\phi\"1\Cb\xin_{k-1}\:\,.
\endgather
$$
Set $\>{\pht\>=\!\tsum_{1\le k<m< n}\!\xi_k\>\xi_m}$. \>Since
\ifMag\vv.15>\fi
$\phi\"2\]=\>\pht\)+\phi\"1\xi_n$, it is a simple exercise to replace
relation \(dim) by
$$
\dim\>\phi\"1\)\Cb\xim_{(\ell-1)}\:+\)
\dim\>\pht\,\Cb\xim_{(\ell-2)}\:\,=\,\dim\)\Si^-\Vlm.
$$
Calculating dimensions explicitly, we further reduce the last relation to
the following one:
\ifMag
$$
\gather
\\
\cnn-1.8>
\dim\>\pht\,\Cb\xim{\vp]}_{(\ell-2)}\:\,=\,\biggl\lb\)
\matrix
\dim\>\Cb\xim{\vp]}_{(\ell-2)}\:\ &\text{if} &\! 2\ell\le n\\
\nn1>
\dim\>\Cb\xim{\vp]}_\ell\:\hfill &\text{if} &\! 2\ell>n
\endmatrix\;,
\\
\nn2>
\Text{which another look is}
\nn6>
{\align
\dim\> & \pht\,\Cb\xim{\vp]}_{(\ell-2)}\:\,={}
\Tag{dim2}
\\
\ifMag\nn4>\else\nn6>\fi
&{}=\,\min\,\bigl(\}\dim\>\Cb\xim{\vp]}_{(\ell-2)}\:\,,\,
\dim\>\Cb\xim{\vp]}_\ell\:\)\bigr)\,.
\kern-2em
\endalign}
\endgather
$$
\else
$$
\align
\dim\>\pht\,\Cb\xim{\vp]}_{(\ell-2)}\:\,=\,{} & \biggl\lb\)
\matrix
\dim\>\Cb\xim{\vp]}_{(\ell-2)}\:\ &\text{if} &\! 2\ell\le n\\
\nn1>
\dim\>\Cb\xim{\vp]}_\ell\:\hfill &\text{if} &\! 2\ell>n
\endmatrix\;,
\\
\Text{which another look is}
\nn6>
\dim\>\pht\,\Cb\xim{\vp]}_{(\ell-2)}\:\,=\,\min &\,\bigl(\}
\dim\>\Cb\xim{\vp]}_{(\ell-2)}\:\,,\,\dim\>\Cb\xim{\vp]}_\ell\:\)\bigr)\,.
\kern-2.4em
\Tag{dim2}
\endalign
$$
\fi
\par
Now it is convenient to change the \var/s. Set
$$
\zt_k=\xi_{k}\lsym+\xi_{n-k-1}\,,\qqq \zt_{n-k}=\xi_{k+1}\lsym+\xi_{n-k}\,,
\qqq 1\le k<n/2 \,,
$$
and $\zt_{n/2}=\xi_{n/2}$ if $n$ is even. Then
$\pht\>=\!\sum_{1\le k<n/2} \zt_k\>\zt_{n-k}$.
\par
Denote by $\zth_m$ and $\phh$ the operators of multiplication by $\zt_m$ and
$\pht$, \resp/. Introduce the Grassman derivations $\der_1\lc\der_{n-1}$
corresponding to $\zt_1\lc\zt_{n-1}$:
$$
\der_k\)1\>=\>0\,,\qqq\der_k\>\zth_m+\)\zth_m\>\der_k\>=\>\dl_{km}\,,\qqq
k,\mn-1\,.
$$
Set $D\>=\!\sum_{1\le k<n/2} \der_k\>\der_{n-k}$. Then the operators
$D\),\,\phh$ and $[D\),\)\phh\>]$ constitute an $\gsl\}$ action in $\Cb\xim$,
and the relation \(dim2) follows from the representation theory of $\gsl$.
\vvgood
\qed
\Pf of formulae \(Xp). Set
\ifMag\vv-.5>\else\vv->\fi
$$
\Tht_k(t)\,=\,\prod_{m=1}^k\,{\Exp(t-z_m)+1\over\Exp(t-z_m)-1}\;.
$$
\vvv.3>
The first of the formulae \(Xp) immeadiately results from the equality
$\Tht_k=\)\Tht_{k-1}+2\)W\'k$. The second formula can be proved by induction
\wrt/ $n$ using the same relation. The base of induction is $n=1$. In this case
the second formula reads \,$\Xp2\!=0$, and is verified straightforwardly.
\epf

\Sect[Uq]{Kernel of the \hmap/ at $\mu=0$ and $\Uq$ at $q=e^{\pi i/2}$}
Consider the \qg/ $\Uq$ with generators $e,\)f,\>k$ subject to relations
$$
\gather
\\
\cnn-1.5>
k\)e\>=\>q^2\)e\)k\,,\qqq kf\>=\>q^{-2}fk\,,\qqq
[\)e\>,\]f\)]\,=\;{k-k\1\over q-q\1\!}\;,
\\
\nn-2>
\Text{and a coproduct $\Dl$:}
\nn5>
\Dl(e)\>=\>1\ox e+e\ox k\,,\qqq \Dl(f)\>=\>k\1\}\ox f+f\ox1\,,\qqq
\Dl(k)\>=\>k\ox k\,,
\\
\endgather
$$
\par
Make the space $V\!$ into a \$\Uq$-module, the generators $e,\)f,\>k$ acting as
\vv.1>
$\si^+\!,\,\si^-\!,\,q^{\si^3}\!\!$, \resp/. Then the generator $f$ is
represented in the \$\Uq$-module $\Von\!$ by the operator $F(q)\in\End(\Von)$:
$$
F(q)\,=\tsun\>q^{-\si^3_1}\}\ldots\>q^{-\si^3_{m-1}}\,\si^-_m\,.
$$
Looking at the limit $q\to e^{\pi i/2}\!$ define operators
$\Fp1\},\)\Fp2\in\End(\Von)$:
$$
\gather
\\
\cnn-1.4>
\Fp1=\)-\)i\>F(e^{\pi i/2})\,=
\tsun(-i)^m\>\si^3_1\ldots\si^3_{m-1}\,\si^-_m\,,\qqq
\\
\nn4>
\Fp2=\,-\)\limq\;{F(q)\vpb2\over 1+q^2}\;=\tsum_{1\le k<m\le n}\!
(-i)^{k+m}\>\si^-_k\)\si^3_{k+1}\ldots\si^3_{m-1}\>\si^-_m\,,
\endgather
$$
the overall normalization factors being fixed just for the aesthetic reason.
There is a faithful \rep/ of $\Cb\xin$ in $\Von\!$ given by
$$
\xi_m\>\map\>(-i)^m\>\si^3_1\ldots\si^3_{m-1}\,\si^-_m\,,\qqq\mn\,,\vp{\Big|}
\Tag{JW}
$$
the operators $\Fp1\},\)\Fp2\!$ representing the elements
$\phi\"1\},\)\phi\"2\!\in\Cb\xin$, \cf. \(phi12). This \rep/ is \eqv/ to
the left regular \rep/ of $\Cb\xin$ with the following intertwiner:
\ifMag\vvn->\fi
$$
\xi_{m_1}\ldots\xi_{m_k}\,\map\,(-i){\vp{\bigg|}}^{\sum_{a=1}^k m_a}v_M\,,
\qqq 1\le m_1\lsym<m_k\le n\,,\kern-2em
\Tag{xiM}
$$
\vvv-.3>
where ${M=\lb\)m_1\lc m_k\)\rb}$ and the vector $v_M$ is defined by \(vM).
In fact, the \rep/ \(JW) is the Jordan-Wigner transformation.
\vsk.3>
For generic $q$ it is known that \>${\dim F(q)\>\Vk=\}\dim\)\Si^-\Vk}$ for any
$k=0\lc n$.  Consider a subspace \>$\limF$, taking the limit in the topo\-logy
of the corresponding Grassmanian. The limit exists because $F(q)$ is \hol/
at $q=i$.
\par
The subspace \>$\limF$ contains both $\Fp1\Vlm$ and\ifMag\nl\fi\ $\Fp2\Vlmm$.
Hence, the equality \(dim) implies that
$$
\limF\,=\,\Fp1\Vlm+\Fp2\Vlmm\,.
\Tag{FFF}
$$
Denote by $C_z$ the composition of \iso/s \(BF) and \(xiM):
\vvn-.5>
$$
C_z:\)\Plus_{k=0}^n\>\Fqbw k(z)\,\to\,\Von\,,\qquad\qquad
C_z:\)W_M(\cdot\>,z)\,\map\,(-i){\vp{\bigg|}}^{\sum_{a=1}^k m_a}v_M\,,
$$
Under this map the operator $X\"a_{z,\ell}$ translates modulo proportionality
into the operator $\Fp a\!$ acting on $(\Von)_{(\ell-a)}\:$. So, by Theorem
\[mu=00] and relation \(FFF) the kernel of the \hmap/ $\Ib_{z,\ell}$
corresponds to the limit of the subspace $F(q)\>\Vlm$ as $q$ tends to $i$.
\vsk.3>
Appearance of the map $C_z$ in the descriptiopn of \sol/s of the \qKZe/ is not
accidental. Utilizing the general picture developed in \Cite{TV}, \Cite{MV} one
can see that for generic step $p$ the space of \sol/s of the \qKZe/ \(qKZ) at
${\mu\ne 0}$ with values in $\Von$ is naturally isomorphic to a tensor product
${\FF\oxc\Von}\!$ of the field of quasiconstants $\FF$ and the \$\Uq$-module
$\Von$ desribed at the beginning of this section, $q$ being equal to
$e^{\pi i\)\h/p}$. The \iso/ is based on the construction of \hsol/s of
the \qKZe/ and its realization in the case ${p=2\h}$, studied here,
is the map $C_z$.
\par
{\bls=1.1\bls
At the same time, for generic step $p$, the space of \sol/s of the \qKZe/
\(qKZ) at ${\mu=0}$ with values in $\Vns$ is isomorphic to either
${\FF\oxc\Vqs}\!$, where $\Vqs\!$ is a subspace of singular vectors \wrt/
the \qg/ $\Uq$, or ${\FF\oxc\bigl(\Von\}\big/F(q)\)\Von\bigr)}$.
These models of the space of \sol/s are \eqv/, since
$\Von=\)\Vqs\oplus\)F(q)\)\Von$ for generic $q$.
\goodbm
\vsk.16>}
The results of this paper means that for $p=2\h$, so that $q=e^{\pi i/2}\!$,
the second model survives with the minimal required modification; namely,
the space of \sol/s of the \qKZe/ \(qKZ) at $\mu=0$ is isomorphic to a tensor
product of $\FF$ and the quotient space
\ifMag\vvn.2>\fi
$$
\Von\}\big/\limq\](F(q)\>\Von\bigr)\,=\,
\Von\}\big/\bigl(\Fp1\Von+\Fp2\Von\bigr)\,.
$$
\ifMag\vvv.3>\fi
Unfortunately, the first model does not survive at $p=2\h$, at least for even
$n$, because in this case the subspaces \>$\limq\]\Vqs\]$ and
\>$\limq\](F(q)\>\Von\bigr)$ intersect nontrivially, and therefore, their sum
does not coincide with $\Von\!$.

\Appendix

\sect{Appendix}
\vsk-.5>\vsk0>
\Pf of Theorem \[dermu].
Let $D$ be the operator defined by \(Df) and denote by $D_a$ the operator $D$
acting \wrt/ the \var/ $t_a$. The operators $D_1\lc D_n$ naturally transform
under \perm/s of the \var/s $\tell$; in particular,
$$
\Asym\bigl(D_1f(\tell)\bigr)\,=
\susi\sgn(\si)\>D_{\si_1}\]\bigl(f(t_{\si_1}\lc t_{\si_\ell})\bigr)\,.
$$
\par
Take a subset $N=\lb\)n_1\lsym<n_{\ell-1}\)\rb$ and consider a total \dif/
$$
\h\1\!\Asym\bigl(D_1\bigl(\)g_N(\twll)\)\tprod_{a=2}^\ell(t_1-t_a-\h)
\bigr)\bigr)\,,
$$
\cf. \(wM). Set ${n_0=0}$. Given $b$ and $m$ \st/ ${1\le b<\ell}$,
${n_{b-1}<m<n_b}$, using Lemma \[D1] we transform this \fn/ as follows:
\ifMag
$$
\align
\h\1(e^\mu\}-1)\>z_m\)w_{N\cup\lb m\rb}\)+\,e^\mu\)w_{N\cup\lb m\rb}\>+
\!\!\sum_{\tsize{k\nin N\atop 1\le k<m}}^{}\!\!\!w_{N\cup\lb k\rb}\>+\,
e^\mu\!\!\!\sum_{\tsize{k\nin N\atop m<k\le n}}\!\!\!w_{N\cup\lb k\rb}\>-{}
\\
\nn9>
{}-\,\h\1(e^\mu\}-1)\>\Asym\}(t_b\,g_{N\cup\lb m\rb})\, &.
\\
\cnn.2>
\endalign
$$
\else
$$
\h\1(e^\mu\}-1)\>z_m\)w_{N\cup\lb m\rb}\)+\,e^\mu\)w_{N\cup\lb m\rb}\>+
\!\!\sum_{\tsize{k\nin N\atop 1\le k<m}}\!\!\!w_{N\cup\lb k\rb}\>+\,
e^\mu\!\!\!\sum_{\tsize{k\nin N\atop m<k\le n}}\!\!\!w_{N\cup\lb k\rb}
-\,\h\1(e^\mu\}-1)\>\Asym\}(t_b\,g_{N\cup\lb m\rb})\,.
$$
\fi
For a subset $M=\lb\)m_1\lsym<m_\ell\)\rb$ let us substitute
$N=M\!\setminus\]\!\lb\)m_a\)\rb$, $b=a$, $m=m_a$ into the last expression
and then take a sum of the results for all $\aell$. The sum equals
\ifMag
$$
\align
\\
\cnn-1.5>
\h\1(e^\mu\}-1)\tsum_{a=1}^\ell\}z_{m_a}\)w_M\>+\,e^\mu\)\ell\>w_M\>+
\sum_{\%{k\nin M}_{\ssize\%{m\in M}_{\ssize k<m\vru1.2ex>}}}
w_{M\cup\lb k\rb\setminus\lb m\rb}\>+\,e^\mu
\sum_{\%{k\nin M}_{\ssize\%{m\in M}_{\ssize k>m\vru1.2ex>}}}
w_{M\cup\lb k\rb\setminus\lb m\rb}\)-{} &
\\
\nn-2>
{}-\,\h\1(e^\mu\}-1)\tsum_{a=1}^\ell\}t_a\,w_M\, &.
\\
\cnn-.5>
\endalign
$$
\else
$$
\gather
\h\1(e^\mu\}-1)\tsum_{a=1}^\ell\}z_{m_a}\)w_M\>+\,e^\mu\)\ell\>w_M\>+
\sum_{\%{k\nin M}_{\ssize\%{m\in M}_{\ssize k<m\vru1.2ex>}}}
w_{M\cup\lb k\rb\setminus\lb m\rb}\>+\,e^\mu
\sum_{\%{k\nin M}_{\ssize\%{m\in M}_{\ssize k>m\vru1.2ex>}}}
w_{M\cup\lb k\rb\setminus\lb m\rb}\)-\,
\h\1(e^\mu\}-1)\tsum_{a=1}^\ell\}t_a\,w_M\,.
\\
\cnn-.7>
\endgather
$$
\fi
We denote it by $r_M$.
\goodbreak
\par
The \difl/ \eq/ \(mueq) for the \hsol/ $\Psi_W$, \cf. \(PsiW) of the \qKZe/ can
be written in the following form
\ifMag\vvn.15>\fi
$$
\Il(r_M,W)\,=\,0\qqq\text{for any}\ \ M\), \ \#M=\ell\,.\kern-2em
$$
\ifMag\vvv.15>\fi
The \fn/ $r_M$ is a total \dif/ by construction, so the claim follows from
Proposition \[ID].
\mmgood
\epf
\Lm{D1}
Let ${N=\lb\)n_1\lsym<n_{\ell-1}\)\rb}$ be a subset of ${\lb\)1\lc n\)\rb}$.
Set ${n_0=0}$, ${n_\ell=n}$. For any $b, m$ \st/ ${1\le b<\ell}$,
${n_{b-1}<m\le n_b}$ the following relations hold:
$$
\gather
\\
\ifMag\cnn-1.5>\else\cnn-1.3>\fi
\alignedat2
&\h\sum_{\tsize{k\nin N\atop 1\le k<m}}\!w_{N\cup\lb k\rb}(\tell)\,={} &&
\\
\nn-1>
&\quad{}=\,\Asym\biggl(\!\biggl(\;\prod_{1<a\le b}(t_1-t_a-\h)\,-
\prod_{1<a\le b}(t_1-t_a+\h)\prod_{1\le k<m}\!{t_1-z_m-\h\over t_1-z_m}\>
\biggr)\,\x{} &&
\\
\nn3>
&& \Llap{{}\x\prod_{b<a\le\ell}(t_1-t_a-\h)\;g_N(\twll)\biggr)\,} &.
\endalignedat
\\
\ald
\nn-6>
\alignedat2
&\h\sum_{\tsize{k\nin N\atop m\le k\le n}}\!w_{N\cup\lb k\rb}(\tell)\,={} &&
\\
\nn-1>
&\quad {}=\,\Asym\biggl(\!\biggl(\;\prod_{b<a\le\ell}(t_1-t_a-\h)\,-
\prod_{b<a\le\ell}(t_1-t_a+\h)\prod_{m\le k\le n}\!{t_1-z_m-\h\over t_1-z_m}\>
\biggr)\,\x{} &&
\\
\nn8>
&&\Llap{{}\x\prod_{1<a\le b}(t_1-t_a+\h)
\prod_{1\le k<m}\!{t_1-z_m-\h\over t_1-z_m}\;g_N(\twll)\biggr)\,} &.
\endalignedat
\endgather
$$
\endpro
\nt
The statement is a simple generalization of Lemma 3.5 in \Cite{NPT}\).
\Pf of formula \(LK).
\ifMag\vvn.2>\fi
Let $\Kt_1\lc\Kt_n$ be products of the \qKZo/s:
$$
\Kt(\zn)\,=\,K_m(z_1+p\lc z_{m-1}+p, z_m\lc z_n)\ldots K_1(\zn)\,.
$$
\ifMag\vvv.2>\fi
Then the compatibility conditions \(LK) can be written as follows:
\ifMag\vvn.2>\fi
$$
[L(\zn)\>{,}\>\Kt_m(\zn)]\,=\,0\,,\qqq\mn\,.
\Tag{LKt}
$$
\ifMag\vvv.2>\fi
The \YB/ for the \Rm/ $R(x)$ implies that
$$
\gather
\\
\ifMag
\cnn-1.2>
{\alignedat2
\Kt_m(\zn)\> &
\Rlap{R_{10}(z_1-u)\ldots R_{n0}(z_n-u)\>\exp\}(-\mu\)H_0)\,={}} &&
\\
\nn7>
{}=\,{} & R_{m+1,0}(z_{m+1}-u)\ldots{}&& R_{n0}(z_n-u)\>\exp\}(-\mu\)H_0)\,\x{}
\\
\nn6>
&& \Llap{{}\x{}\,R_{10}(z_1-u)\ldots{}} & R_{m0}(z_m-u)\>\Kt_m(\zn)\,,
\endalignedat}
\else
\cnn-1.5>
{\align
\Kt_m(\zn)\>R_{10}(z_1-u)\ldots R_{n0}(z_n-u)\> &\exp\}(-\mu\)H_0)\,={}
\\
\nn6>
{}=\,R_{m+1,0}(z_{m+1}-u)\ldots R_{n0}(z_n-u)\> &\exp\}(-\mu\)H_0)
R_{10}(z_1-u)\ldots R_{m0}(z_m-u)\>\Kt_m(\zn)\,,
\endalign}
\fi
\\
\ifMag\nn-4>\else\nn4>\fi
\Text{therefore,}
\nn2>
\bigl[\Kt_m(\zn)\>,\tr_{\)0}\:\bigl(R_{10}(z_1-u)\ldots R_{n0}(z_n-u)\>
\exp\}(-\mu\)H_0)\bigr)\bigr]\,=\,0\,.
\Tag{Ktr}
\endgather
$$
\ifMag\vvv.1>\fi
On the other hand the operator $L$ appears in the following expansion:
\ifMag
$$
\align
e^\mu\)\tr_{\)0}\:\bigl(R_{10}(z_1-u)\ldots R_{n0}(z_n-u)\>
\exp\}(-\mu\)H_0)\bigr)\,=\,
c_{00}\:+u\1\)\bigl(\)c_{10}\:+c_{11}\:\tsun H_m\bigr)\,+\!{} &
\\
\nn4>
{}+\,u^{-2}\)\bigl(\)c_{20}\:+c_{21}\:\tsun H_m+
c_{22}\:\>\bigl(\tsun H_m\bigr)^2+\h\>(e^\mu\}-1)\,L\)\bigr)+O(u^{-3}) &
\endalign
$$
\else
$$
\align
& e^\mu\)\tr_{\)0}\:\bigl(R_{10}(z_1-u)\ldots R_{n0}(z_n-u)\>
\exp\}(-\mu\)H_0)\bigr)\,={}
\\
\nn6>
&\!\}{}=\,c_{00}\:+u\1\)\bigl(\)c_{10}\:+c_{11}\:\tsun H_m\bigr)+
u^{-2}\)\bigl(\)c_{20}\:+c_{21}\:\tsun H_m+c_{22}\:\>\bigl(\tsun H_m\bigr)^2
+\h\>(e^\mu\}-1)\,L\)\bigr)+O(u^{-3})
\endalign
$$
\fi
where $c_{ij}\:$ are \$\C\)$-valued coefficients depending on $\mu\),\,\zn$.
\vvn.1>
Hence, \(Ktr) entails the required relation \(LKt).
\vvgood\mmgood
\epf
\Pf of inclusion \>${\im X\"2_{z,\ell}\sub\ker\Ib_{z,\ell}}$.
We give here only a very short draft of the proof to point out the technical
piece of subtlety in this case compared with Lemma 5.3 in \Cite{NPT}.
\vsk.3>
For brevity we consider the case $\ell=2$. Set
$$
\gather
\\
\cnn-1.4>
E(t_1,t_2)\,=\;{\Exp(t_1-t_2)-1\over\Exp(t_1-t_2)+1}\;,\qqq
F(t_1,t_2)\,=\,\Tht(t_1)\>\Tht(t_2)\>E(t_1,t_2)\,,
\\
\ifMag\nn2>\fi
\Text{so that}
\ifMag\nn-1>\else\nn-3>\fi
\Xp2(t_1,t_2)\,=\,F(t_1,t_2)-E(t_1,t_2)+\Tht(t_1)-\Tht(t_2)\,.
\endgather
$$
It turns out that for each term $W$ in \rhs/ of the above formula all
the \hint/s $\Itw(w_M,W)$, $\#M=2$. For the third and fourth terms
the proof of this assertion is like in \Cite{NPT}.
\par
For the terms $E$ and $F$ we first have to change the integration \var/s:
${t_1=u_1}$, $t_2=u_1+u_2$. After that, considering the integration \wrt/
the \var/ $u_1$ we can apply arguments like in \Cite{NPT} to prove vanishing
of the \hint/s $\Itw(w_M,E)$ and $\Itw(w_M,F)$. Notice that the integration
contour in the definition of these integarls, \cf. \(hint), can be chosen
avoiding extra singularities of the integrands produced by the term
$\Exp(t_1-t_2)+1$ in the denominators.
\epf

\myRefs
\widest{VW}
\parskip.1\bls

\ref\Key FTV
\by \Feld/, \VT/ and \Varch/ \info in preparation
\endref

\ref\Key MV
\by E\&Mukhin and \Varch/
\paper The quantized \KZv/ \eq/ in tensor products of \irr/ \$sl_2$-modules
\jour Preprint \yr 1997 \pages 1--32
\endref

\ref\Key NPT
\by A\&Nakayashiki, S\&Pakuliak and \VT/
\paper On \sol/s of the \KZ/ and \qKZ/ equations at level zero
\jour Preprint \yr 1997 \pages 1--24
\endref

\ref\Key R
\by \Reshy/
\paper Jackson-type integrals, \Bv/s, and \sol/s to a \dif/ analogue
of the \KZv/ system
\jour \LMP/ \vol 26 \yr 1992 \pages 153--165
\endref

\ref\Key S1
\by F\]\&A\&Smirnov
\book Form factors in completely integrable models of quantum field theory
\bookinfo Advanced Series in Math.\ Phys., vol\&14
\yr 1992 \publ \WSa/
\endref

\ref\Key S2
\by F\]\&A\&Smirnov
\paper On the deformation of Abelian integrals
\jour \LMP/ \vol 36 \yr 1996 \pages 267--275
\endref

\ref
\by\refin F\]\&A\&Smirnov
\paper Counting the local fields in SG theory
\jour Nucl.\ Phys. \vol B453 \yr 1995 \pages 807--824
\endref

\ref\Key TV
\by \VT/ and \Varch/
\paper Geometry of \$q$-hypergeometric \fn/s as a bridge between Yangians and
\qaff/s \jour \Inv/ \yr 1997 \vol 128 \issue 3 \pages 501--588
\endref

\ref\Key V
\by \Varch/
\paper Quantized \KZv/ \eq/s, quantum \YB/, and \deq/s for \$q$-\hgeom/ \fn/s
\jour \CMP/ \vol 162 \yr 1994 \pages 499--528
\endref

\ref\Key WW
\by E\&T\]\&Whittaker and G\&N\&Watson
\book A Course of Modern Analysis
\yr 1927 \publ \CUP/
\endref

\endRefs

\bye